\documentclass[11pt,a4paper]{amsart}

\usepackage[T1]{fontenc}
\usepackage[utf8]{inputenc}
\usepackage{lmodern}
\usepackage{microtype}
\usepackage{amsmath,amssymb,amsthm,mathtools}
\usepackage{mathrsfs}
\usepackage{enumitem}
\usepackage{etoolbox}
\patchcmd{\part}{\raggedright}{\centering}{}{%
  \PackageError{integral-bc}{Could not center part headings}
    {Check the amsart definition of \string\part.}}

\usepackage{booktabs}
\usepackage{array}
\usepackage{xcolor}
\usepackage[a4paper,margin=31mm]{geometry}
\usepackage[
  colorlinks=true,
  linkcolor=blue!55!black,
  citecolor=green!35!black,
  urlcolor=blue!60!black,
  pdfauthor={Bora Yalkinoglu},
  pdftitle={Integral and p-adic Bost-Connes systems}
]{hyperref}

\newcommand{\A}{\mathcal A}
\newcommand{\C}{\mathbb C}

\newcommand{\F}{\mathbb F}

\newcommand{\Kbar}{\overline{K}}
\newcommand{\N}{\mathbb N}
\newcommand{\Ocal}{\mathcal O}
\newcommand{\Ok}{\mathcal O_K}
\newcommand{\Obar}{\overline{\mathcal O}_K}
\newcommand{\Q}{\mathbb Q}

\newcommand{\W}{\mathbb W}
\newcommand{\Z}{\mathbb Z}
\newcommand{\ab}{\mathrm{ab}}
\newcommand{\DR}{\mathrm{DR}}
\newcommand{\Gal}{\operatorname{Gal}}
\newcommand{\Hom}{\operatorname{Hom}}

\newcommand{\Clc}{C_{\mathrm{lc}}}

\newcommand{\red}{\mathrm{red}}
\newcommand{\et}{\mathrm{\acute{e}tale}}
\newcommand{\id}{\mathrm{id}}

\newcommand{\Nm}{\mathrm N}
\newcommand{\tensor}{\mathbin{\otimes}}
\newcommand{\into}{\lhook\joinrel\longrightarrow}
\newcommand{\onto}{\relbar\joinrel\twoheadrightarrow}

\newtheorem{theorem}{Theorem}[section]
\newtheorem{proposition}[theorem]{Proposition}
\newtheorem{corollary}[theorem]{Corollary}
\newtheorem{lemma}[theorem]{Lemma}
\theoremstyle{definition}
\newtheorem{definition}[theorem]{Definition}
\newtheorem{example}[theorem]{Example}
\theoremstyle{remark}
\newtheorem{remark}[theorem]{Remark}

\numberwithin{equation}{section}

\title[Integral and $p$-adic Bost-Connes systems]
{On Integral and $p$-adic Bost-Connes systems}

\author{Bora Yalkinoglu}
\address{CNRS and IRMA, Universit\'e de Strasbourg, 7 rue Ren\'e Descartes,
67084 Strasbourg, France}
\email{yalkinoglu@math.unistra.fr}
\date{\today}

\subjclass[2020]{11R23, 11R37, 11R42, 13F35, 46L55}
\keywords{Bost-Connes system, Deligne-Ribet monoid,
periodic Witt vectors, Verschiebung, integral model, \(p\)-adic \(L\)-function, \(p\)-adic KMS condition}

\begin{document}

\begin{abstract}
Using Borger and de Smit's theory of periodic Witt vectors, we construct integral models of Bost-Connes
systems for general number fields, generalizing Connes, Consani and
Marcolli's construction over \(\Q\).
As an application, for totally real number fields other than \(\Q\), we construct $p$-adic Bost-Connes systems whose \(p\)-adic KMS states are described in terms of Deligne-Ribet's theory of $p$-adic $L$-functions, extending the work of Connes and Consani for $\mathbb Q$.
\end{abstract}

\maketitle

\section{Introduction}
\label{sec:introduction}

\noindent Bost-Connes systems were introduced in the seminal paper
\cite{BC95}. For subsequent developments see, e.g.,
\cite{HP05,CMR05,LLN09,CCM09,Yal13,CC14}.
To each number field \(K\) one can associate a Bost-Connes
system \(\A_K=(A_K,\sigma_t)\), a \(C^*\)-dynamical system
with the following key properties:
\begin{enumerate}[label=(\roman*),leftmargin=*,itemsep=3pt]
\item Its partition function is the Dedekind zeta function
\(\zeta_K\).
\item The maximal abelian Galois group \(\Gal(K^{\ab}/K)\)
provides symmetries of \(\A_K\).
\item There is a \(K\)-rational subalgebra
\(A_K^{\mathrm{arith}}\subset A_K\) such that, for every extremal
\(\mathrm{KMS}_\infty\)-state \(\varrho\), the values
\(\varrho(f)\), \(f\in A_K^{\mathrm{arith}}\), belong to
\(K^{\ab}\) and generate it over \(K\).
\item For \(\nu\in\Gal(K^{\ab}/K)\) and $f \in A^{\mathrm{arith}}_K$, the following compatibility holds:
\[
 {}^\nu\varrho(f)=\nu^{-1}\bigl(\varrho(f)\bigr).
\]
\item The \(\C\)-span of \(A_K^{\mathrm{arith}}\) is dense in \(A_K\). \\ 

\end{enumerate} 
\noindent We recall the construction of \(A_K^{\mathrm{arith}}\)
from \cite{Yal13}, which is based on the Grothendieck-Galois correspondence of \(\Lambda_K\)-rings from \cite{BdS08,BdS18} and the endomotive formalism of \cite{CCM07}. Details are given in Section~\ref{sec:DR}. \\ \\
\noindent Let \(\Ok\) be the ring of integers of \(K\), let
\(\Kbar\) be an algebraic closure, and let \(\Obar\) be the
integral closure of \(\Ok\) in \(\Kbar\).
Write \(I_K\) for the monoid of nonzero integral ideals and
\(K(\mathfrak f)\subset K^{\mathrm{ab}}\) for the strict ray class field of modulus
\(\mathfrak f\), with all real places included.
Under the Grothendieck-Galois correspondence for $\Lambda_K$-rings, each finite Deligne-Ribet monoid
\(\DR_{\mathfrak f}\), \(\mathfrak f\in I_K\), corresponds to a finite \'etale \(K\)-algebra \(E_{\mathfrak f}\). Explicitly,
\[\DR_{\mathfrak f} \cong \Hom_{K\text{-alg}}(E_{\mathfrak f},\Kbar)
 ,\qquad
 E_{\mathfrak f}\cong
 \prod_{\mathfrak d\mid\mathfrak f}
       K(\mathfrak f\mathfrak d^{-1}).
\]
In the limit, we have the profinite Deligne-Ribet monoid $\DR_K$ and the commutative part of $A^{\mathrm{arith}}_K$, the $K$-algebra $E_K$ defined by
\[
 \DR_K:=\varprojlim_{\mathfrak f}\DR_{\mathfrak f}
 \cong\Hom_{K\text{-alg}}(E_K,\Kbar), \qquad  E_K:=\varinjlim_{\mathfrak f}E_{\mathfrak f}.
\]

\noindent Multiplication by \(\mathfrak a\in I_K\) defines the clopen subset \(\mathfrak a\DR_K \subset \DR_K\), which corresponds to the
idempotent \(e_{\mathfrak a}\in E_K\).
It induces the unital Frobenius endomorphism
\(\sigma_{\mathfrak a}\) of \(E_K\), whose restriction to
\(e_{\mathfrak a}E_K\) is an isomorphism.  Define
\[
 \rho_{\mathfrak a}
 :=\bigl(\sigma_{\mathfrak a}|_{e_{\mathfrak a}E_K}\bigr)^{-1}
 :E_K\longrightarrow e_{\mathfrak a}E_K\subseteq E_K.
\]
Thus \(\sigma_{\mathfrak a}\rho_{\mathfrak a}=\id_{E_K}\) and
\(\rho_{\mathfrak a}(1)=e_{\mathfrak a}\). \\ \\
\noindent The maps $\rho_{\mathfrak a}$ are then used to define the arithmetic subalgebra $A_K^{\mathrm{arith}}$ given by the crossed product
\[
 A_K^{\mathrm{arith}}=E_K\rtimes_\rho I_K.
\]

\noindent The goal\footnote{Our results realize the integral model $A_K^{\mathrm{int}}$ and the connection to Deligne-Ribet's $p$-adic \(L\)-functions as promised in \cite{Yal13}.} of the present paper is to construct integral $A_K^{\mathrm{int}}$ and \(p\)-adic $A_{K,p}$ versions of \(A_K^{\mathrm{arith}}\),
continuing, e.g., the works \cite{CC14,CC21,CC22}. \\ \\
\noindent Our paper has two parts. In Part~\ref{part:integral} we use the theory of periodic Witt vectors 
of Borger and de Smit
\[
 R_{\mathfrak f}:=\W_K^{(\mathfrak f)}(\Obar) \subset E_\mathfrak f
\]
to define an $\Ok$-algebra, see Section~\ref{sec:orders}, \[R_K:=\varinjlim_{\mathfrak f}R_{\mathfrak f} \subset E_K.\]
The maps \(\sigma_{\mathfrak a}\) naturally preserve \(R_K\), but the \(\rho_{\mathfrak a}\) do not, i.e., we cannot immediately define an integral version of $A_K^{\mathrm{arith}}$ using the maps $\rho_{\mathfrak a}$.

\noindent This obstruction already occurs over \(\Q\), where Connes,
Consani and Marcolli \cite[\S4]{CCM09} replace \(\rho_n\) by the integral additive operator \(n\rho_n\), $n \in \N$. \\ \\
\noindent For general \(K\), the replacement is the family
\[
 V_{\mathfrak a,\xi}(x)=\xi\rho_{\mathfrak a}(x),
 \qquad \xi\in\mathfrak a,\quad x\in R_K.
\]
We prove that these scaled maps  preserve \(R_K\): they are
the restrictions of Borger's \(\Ok\)-linear Verschiebung operators
on \(\W_K(\Obar)\) \cite[\S\S4.2--4.3]{Borger11}.
The family avoids choosing a generator of a possibly nonprincipal
ideal.  For \(K=\Q\), \(\mathfrak a=(n)\) and \(\xi=n>0\),
it recovers the usual \(V_n\), whose restriction to \(R_{\Q}\)
is \(n\rho_n\).

\noindent Write \(U_{\mathfrak a},U_{\mathfrak a}^*\) for the usual
generators of \(A_K^{\mathrm{arith}}\), with
\(U_{\mathfrak a}^*U_{\mathfrak a}=1\) and
\(\rho_{\mathfrak a}(x)=U_{\mathfrak a}xU_{\mathfrak a}^*\). \\ \\
\noindent We define \(A_K^{\mathrm{int}}\) to be the \(\Ok\)-subalgebra
generated by \(R_K\), \(U_{\mathfrak a}^*\), and
\(\xi U_{\mathfrak a}\) for \(\xi\in\mathfrak a\).
Here \(\xi U_{\mathfrak a}\) is simply scalar multiplication
inside the known arithmetic algebra. The generator \(U_{\mathfrak a}\) alone
need not belong to the integral algebra.

\begin{theorem}
\label{thm:intro-integral}
\noindent The algebra \(A_K^{\mathrm{int}}\) has the following properties.
\begin{enumerate}[label=(\roman*),leftmargin=*,itemsep=3pt]
\item It is flat over \(\Ok\), and
\[
 A_K^{\mathrm{int}}\cap E_K=R_K,\qquad
 K\tensor_{\Ok}A_K^{\mathrm{int}}\cong A_K^{\mathrm{arith}}.
\]
\item Inside \(A_K^{\mathrm{arith}}\), its normal form is
\[
 A_K^{\mathrm{int}}
 =\bigoplus_{\substack{\mathfrak a,\mathfrak b\in I_K\\
                      (\mathfrak a,\mathfrak b)=1}}
       U_{\mathfrak a}(\mathfrak aR_K)U_{\mathfrak b}^*.
\]
Thus every element has a unique expression as a finite sum
\(\sum U_{\mathfrak a}z_{\mathfrak a,\mathfrak b}U_{\mathfrak b}^*\),
where \(z_{\mathfrak a,\mathfrak b}\in\mathfrak aR_K\)
and the pairs are coprime. Here \(\mathfrak aR_K\) denotes the ideal
of \(R_K\) generated by \(\mathfrak a\).
\item For \(\xi\in\mathfrak a\), Verschiebung restricts to
\[
 V_{\mathfrak a,\xi}:R_{\mathfrak f}\longrightarrow
 R_{\mathfrak a\mathfrak f},\qquad
 V_{\mathfrak a,\xi}(x)=\xi\rho_{\mathfrak a}(x).
\]
For \(K=\Q\), these operators and \(A_K^{\mathrm{int}}\)
recover the integral BC-algebra of Connes, Consani and Marcolli
\cite{CCM09}.
\end{enumerate}
Thus \(A_K^{\mathrm{int}}\) is an integral model of \(A_K^{\mathrm{arith}}\).
\end{theorem}

\noindent 
The theorem follows from Proposition~\ref{prop:V-periodic},
Theorem~\ref{thm:witt-integral-algebra}, and
Corollaries~\ref{cor:concrete-integral} and
\ref{thm:Q-specialization}. \\ \\
\noindent In Part~\ref{part:padic}, we fix a prime number \(p\) and an embedding \(\jmath:\Kbar\into\C_p\).
Write \(I_K^{(p)}=\{\mathfrak a\in I_K:(\mathfrak a,p\Ok)=1\}\),
and let \(A_K^{\mathrm{int},(p)}\subseteq A_K^{\mathrm{int}}\)
be the \(\Ok\)-subalgebra generated by \(R_K\),
\(U_{\mathfrak a}^*\), and \(\xi U_{\mathfrak a}\), for
\(\mathfrak a\in I_K^{(p)}\) and \(\xi\in\mathfrak a\).
The \(p\)-adic BC-algebra is
\begin{equation}
 A_{K,p}:=\C_p\tensor_{\Ok,\jmath}A_K^{\mathrm{int},(p)}.
\label{eq:intro-padic-algebra}
\end{equation}

\noindent The \(p\)-adic BC-system consists of \(A_{K,p}\) and the
automorphism \(\sigma^{(\beta)}\), depending on \(\beta \in D_{K,p}\), the $p$-adic disk defined in section \ref{subsec:padic-dynamics}. They play the role of the $p$-adic time evolution. 
Each automorphism fixes \(\C_p\tensor_{\Ok,\jmath}R_K\) and satisfies
\[
 \sigma^{(\beta)}(U_{\mathfrak a}^*)
       =(\Nm(\mathfrak a))^{(\beta)}U_{\mathfrak a}^*,\qquad
 \sigma^{(\beta)}(\xi U_{\mathfrak a})
       =\bigl((\Nm(\mathfrak a))^{(\beta)}\bigr)^{-1}
          \xi U_{\mathfrak a}.
\]
A \(\C_p\)-linear form \(\varphi:A_{K,p}\to\C_p\) satisfies
the \emph{\(p\)-adic KMS condition at parameter \(\beta\)} if
\begin{equation}
 \varphi\bigl(X\sigma^{(\beta)}(Y)\bigr)=\varphi(YX)
 \qquad(X,Y\in A_{K,p}).
\label{eq:intro-KMS}
\end{equation}
Such a form is a \emph{\(p\)-adic KMS state} when
\(\varphi(1)=1\). No positivity condition is imposed.

\noindent A finite-order totally even ray character \(\chi\) defines an element
\(x_\chi\in\C_p\tensor_K E_K\) 
(Section~\ref{subsec:DR-main}).

\begin{theorem}
\label{thm:intro-padic}
\noindent Let \(K\ne\Q\) be totally real.
For \(\beta\in D_{K,p}\setminus\{1\}\), the Deligne-Ribet measures
define \(p\)-adic KMS forms
\(\varphi_\beta:A_{K,p}\to\C_p\), satisfying
\eqref{eq:intro-KMS}.
For each \(A\in A_{K,p}\), the function
\(\beta\mapsto\varphi_\beta(A)\) is analytic away from \(1\),
with at most a simple pole there.
Their arithmetic values are
\[
 \varphi_\beta(1)=\zeta_{K,p}(\beta),\qquad
 \varphi_\beta(x_\chi)=L_{K,p}(\beta,\chi)
\]
for every finite-order totally even ray character \(\chi\),
with the normalization of
Section~\ref{subsec:DR-main}.
For \(\beta\ne1\) with \(\zeta_{K,p}(\beta)\ne0\), the form
\[
 \widehat\varphi_\beta:=
       \zeta_{K,p}(\beta)^{-1}\varphi_\beta
\]
is a \(p\)-adic KMS state.
If \(K\) is real quadratic, this normalized family extends
pointwise analytically to \(\beta=1\), where it is again a
\(p\)-adic KMS state.
\end{theorem}

\noindent Theorem~\ref{thm:DR-main-KMS} and
Corollary~\ref{cor:quadratic-weight-one} give the construction
and proofs.
\subsection*{Outline}
\noindent The paper is organized as follows.
Section~\ref{sec:DR} recalls our construction of Bost-Connes endomotives from \cite{Yal13}.
Section~\ref{sec:orders} describes periodic Witt vectors and their
Frobenius and Verschiebung operators from \cite{Borger11,BdS18}. Section~\ref{sec:W-algebra}
constructs the integral algebra and its Witt-vector representation.
Section~\ref{sec:padic} then passes to the definition of \(p\)-adic BC-systems and the important $p$-adic KMS condition following \cite{CC14},
and Section~\ref{sec:DR-padic} constructs $p$-adic KMS states using
Deligne-Ribet's theory of $p$-adic measures and $L$-functions for totally real number fields \cite{DR80} and describes their arithmetic properties.

\section*{Notation}

\noindent Throughout, \(K\) is a number field with ring of integers
\(\Ok\), \(\Kbar\) is an algebraic closure, and \(\Obar\) is the
integral closure of \(\Ok\) in \(\Kbar\).
We write \(G_K=\Gal(\Kbar/K)\) and
\(G_K^{\ab}=\Gal(K^{\ab}/K)\).
The monoid \(I_K\) of nonzero integral ideals has group completion
\(J_K=I_K^{\mathrm{gp}}\), the group of nonzero fractional ideals.
The absolute norm \(\Nm(\mathfrak a)=|\Ok/\mathfrak a|\) extends by
\(\Nm(\mathfrak a\mathfrak b^{-1})
=\Nm(\mathfrak a)/\Nm(\mathfrak b)\).
All ray moduli are strict: their infinite part contains every real
place.  We denote their finite part by \(\mathfrak f\), and the
corresponding ray class field by \(K(\mathfrak f)\).

\noindent Write \(\mathbb A_{K,f}\) for the finite adeles,
\(\widehat{\mathcal O}_K=\prod_{\mathfrak q}\Ocal_{K_{\mathfrak q}}\), and
\(\operatorname{Art}_K:\mathbb A_{K,f}^{\times}\to G_K^{\ab}\)
for the Artin map with the convention of \cite{Yal13}.
The letter \(p\) denotes a rational prime, and \(\mathfrak p,\mathfrak q\)
denote prime ideals.  In the \(p\)-adic constructions,
\(\C_p\) is the completed algebraic closure of \(\Q_p\), with
\(|p|_p=p^{-1}\), and an embedding \(\jmath:\Kbar\into\C_p\)
is fixed.

\part{Integral models}
\label{part:integral}

\section{The Bost-Connes endomotive}
\label{sec:DR}

\noindent This section recalls the algebraic model of
\cite[\S\S3,5-7]{Yal13}, using the Galois correspondence of
\cite[Theorem~1.2 and Corollary~5.6]{BdS18}
and the endomotive formalism of \cite{CCM07}.

\subsection{Ray class algebras and the Deligne-Ribet monoid}

\noindent The inclusion \(I_K\subset J_K\) allows integral-ideal
multiplication to be graded by fractional ideals:
every \(\gamma\in J_K\) has a unique expression
\(\gamma=\mathfrak a\mathfrak b^{-1}\) with
\(\mathfrak a,\mathfrak b\in I_K\) coprime, obtained by separating
the positive and negative prime exponents.
For \(\mathfrak f\in I_K\), the Deligne-Ribet monoid is
\(\DR_{\mathfrak f}=I_K/{\sim_{\mathfrak f}}\), where
\begin{equation}
 \mathfrak a\sim_{\mathfrak f}\mathfrak b
 \quad\Longleftrightarrow\quad
 \mathfrak a=x\mathfrak b
 \text{ for some }x\in K_+^\times\cap(1+\mathfrak f\mathfrak b^{-1}).
\label{eq:DR-equivalence}
\end{equation}
Here \(K_+^\times\) denotes the totally positive elements.
Set \(\DR_K:=\varprojlim_{\mathfrak f}\DR_{\mathfrak f}\).
Its unit group is \(G_K^{\ab}\) under Artin reciprocity.
The adelic description of \cite[Proposition~7.2]{Yal13} is
\begin{equation}
 \DR_K\cong
 \widehat{\mathcal O}_K\times_{\widehat{\mathcal O}_K^{\times}}G_K^{\ab},
 \qquad
 u\cdot(r,g)=(ru,\operatorname{Art}_K(u)^{-1}g).
\label{eq:DR-adelic}
\end{equation}

\noindent The finite \(\et\) \(K\)-algebra \(E_{\mathfrak f}\) is
characterized by: 
\begin{equation}
 \Hom_{K\text{-alg}}(E_{\mathfrak f},\Kbar)\cong\DR_{\mathfrak f},
 \qquad
 E_{\mathfrak f}\cong
 \prod_{\mathfrak d\mid\mathfrak f}K(\mathfrak f\mathfrak d^{-1}).
\label{eq:E-f-product}
\end{equation}
These identifications are \(G_K\)-equivariant.  The monoid projections give
injective algebra maps and
\begin{equation}
 E_K:=\varinjlim_{\mathfrak f}E_{\mathfrak f},\qquad
 \DR_K\cong\Hom_{K\text{-alg}}(E_K,\Kbar).
\label{eq:E-K-points}
\end{equation}
In particular, \(\DR_K\) is a profinite topological monoid.
For \(F\in E_K\) and \(x\in \DR_K\), we write \(F(x):=x(F)\).

\begin{example}
\noindent For \(K=\Q\) and modulus \(n\), the map
\((a)\mapsto a\bmod n\) identifies \(\DR_n\) with the
multiplicative monoid \(\Z/n\Z\).
After choosing a primitive \(n\)-th root \(\zeta_n\),
\[
 E_n=\Q[T]/(T^n-1),\qquad
 \Hom_{\Q\text{-alg}}(E_n,\overline\Q)
 =\{T\mapsto\zeta_n^a:a\in\Z/n\Z\}.
\]
The monoid action is precomposition by \(T\mapsto T^m\).
The units of \(\Z/n\Z\) act through cyclotomic Galois automorphisms.
This is the finite-level form of the Hom description above.
\end{example}

\subsection{The arithmetic subalgebra $A_K^{\mathrm{arith}}$}

\noindent The relation
\begin{equation}
 \mathfrak b\sim_{\mathfrak f}\mathfrak c
 \quad\Longleftrightarrow\quad
 \mathfrak a\mathfrak b\sim_{\mathfrak a\mathfrak f}
 \mathfrak a\mathfrak c
\label{eq:multiplication-ray-equivalence}
\end{equation}
gives an injection
\(\DR_{\mathfrak f}\into\DR_{\mathfrak a\mathfrak f}\).
On the limit, multiplication by \(\mathfrak a\) is injective with
clopen image \(\mathfrak a\DR_K\).
Let \(e_{\mathfrak a}\in E_K\) be the corresponding idempotent.
The unital endomorphisms \(\sigma_{\mathfrak a}\) and partial
isomorphisms \(\rho_{\mathfrak a}:E_K\to e_{\mathfrak a}E_K\)
are determined by
\begin{equation}
 (\sigma_{\mathfrak a}F)(x)=F(\mathfrak ax),\qquad
 (\rho_{\mathfrak a}F)(y)=
 \begin{cases}F(x),&y=\mathfrak ax,\\0,&y\notin \mathfrak a\DR_K.
 \end{cases}
\label{eq:rho-functions}
\end{equation}
They satisfy
\(\sigma_{\mathfrak a}\rho_{\mathfrak a}=\id\) and
\(\rho_{\mathfrak a}\sigma_{\mathfrak a}(F)=e_{\mathfrak a}F\).

\noindent The algebra of the Bost-Connes endomotive is
\begin{equation}
 A_K^{\mathrm{arith}}=E_K\rtimes_\rho I_K.
\label{eq:algebraic-endomotive}
\end{equation}
It is generated by \(E_K\) and \(U_{\mathfrak a},U_{\mathfrak a}^*\),
with
\begin{equation}
\begin{gathered}
 U_{\mathfrak a}^*U_{\mathfrak a}=1,\quad
 U_{\mathfrak a}U_{\mathfrak a}^*=e_{\mathfrak a},\quad
 U_{\mathfrak a}U_{\mathfrak b}=U_{\mathfrak a\mathfrak b},\quad
 U_{\mathfrak a}^*U_{\mathfrak b}^*=U_{\mathfrak a\mathfrak b}^*,\\
 U_{\mathfrak a}FU_{\mathfrak a}^*=\rho_{\mathfrak a}(F),\qquad
 U_{\mathfrak a}^*F=\sigma_{\mathfrak a}(F)U_{\mathfrak a}^*,
\end{gathered}
\label{eq:crossed-relations}
\end{equation}
and \(U_{\mathfrak a}U_{\mathfrak b}^*
=U_{\mathfrak b}^*U_{\mathfrak a}\) for coprime ideals.
The star denotes the second algebraic generator.
The direct \(J_K\)-grading is
\begin{equation}
 A_K^{\mathrm{arith}}=\bigoplus_{\gamma\in J_K}A_{K,\gamma}^{\mathrm{arith}},\qquad
 U_{\mathfrak a}FU_{\mathfrak b}^*
 \in A_{K,\mathfrak a\mathfrak b^{-1}}^{\mathrm{arith}}.
\label{eq:JK-grading}
\end{equation}
The identification with the arithmetic subalgebra of the complex
BC-system is \cite[Proposition~6.9 and Theorems~7.1,7.12]{Yal13}.

\subsection{The analytic BC-system}
\label{subsec:complex-system}

\noindent Write
\(C(\DR_K)\) for the \(C^*\)-algebra of continuous complex-valued
functions. The maps
\(\rho_{\mathfrak a}\) of \eqref{eq:rho-functions} act on
\(C(\DR_K)\) as well.

\begin{definition}
\noindent The BC-system attached to \(K\) is
\[
 \A_K=(A_K,\sigma_t),\qquad
 A_K=C(\DR_K)\rtimes_\rho I_K,
\]
where the crossed product
is the $C^*$-algebra generated by \(C(\DR_K)\) and isometries \(U_{\mathfrak a}\) with
\[
 U_{\mathfrak a}U_{\mathfrak b}=U_{\mathfrak a\mathfrak b},\qquad
 U_{\mathfrak a}FU_{\mathfrak a}^*=\rho_{\mathfrak a}(F).
\]
Its time evolution is
\[
 \sigma_t(F)=F,\qquad
 \sigma_t(U_{\mathfrak a})=(\Nm(\mathfrak a))^{it}U_{\mathfrak a}.
\]
\end{definition}

\noindent This is the system of \cite[\S2]{LLN09}, under
\cite[Corollary~7.4]{Yal13}.
Multiplication by \(g\in \DR_K^\times=\Gal(K^{\ab}/K)\) induces
the symmetry
\(\theta_g(F)(x)=F(gx)\), \(\theta_g(U_{\mathfrak a})=U_{\mathfrak a}\).
After choosing \(\jmath_\infty:\Kbar\into\C\), evaluation
\[
 F\longmapsto\bigl[x\longmapsto\jmath_\infty(x(F))\bigr]
\]
embeds \(E_K\) into \(C(\DR_K)\).
Adjoining the isometries and their adjoints realizes
\(A_K^{\mathrm{arith}}\) as a subalgebra of \(A_K\) whose complex
span is dense \cite[Theorem~7.12]{Yal13}.

\begin{remark}
\noindent The KMS classification is also known
\cite[Theorem~2.1]{LLN09}: the state is unique for
\(0<\beta\le1\), while for \(\beta>1\) the Galois group acts
freely and transitively on the extremal states.
We will distinguish these analytic properties from the
\(p\)-adic statements proved in Section~\ref{sec:DR-padic}.
\end{remark}

\section{Periodic Witt vectors}
\label{sec:orders}

\noindent Generalized and periodic Witt vectors were introduced by Borger and de Smit in
\cite{Borger11,BdS18}. We recall the results needed for our construction, and then verify the
Verschiebung compatibility used in the integral BC-algebra.

\subsection{Deligne-Ribet monoids and periodic Witt vectors}

\noindent A \(\Lambda_K\)-structure on a flat \(\Ok\)-algebra \(A\) consists
of commuting \(\Ok\)-algebra endomorphisms \(\psi_{\mathfrak p}\), one for each prime,
such that \(\psi_{\mathfrak p}(x)\equiv x^{\Nm(\mathfrak p)}
\pmod{\mathfrak p A}\).  Products give \(\psi_{\mathfrak a}\) for
\(\mathfrak a\in I_K\).
For flat \(A\), Borger defines \(\W_K(A)\) as the largest
\(\Ok\)-subalgebra of the ghost ring \(A^{I_K}\) stable under
the shifts and on which they are Frobenius lifts
\cite[\S1]{Borger11}.  This characterization explains the
Frobenius-lift structure, not just the underlying coordinates.
The functor extends to arbitrary \(\Ok\)-algebras. For flat \(A\),
the ghost map is injective and
\begin{equation}
 \W_K(A)\into A^{I_K},\qquad
 (\psi_{\mathfrak a}x)_{\mathfrak b}=x_{\mathfrak a\mathfrak b}.
\label{eq:ghost-shift}
\end{equation}
\begin{definition}
\noindent The ring of \(\mathfrak f\)-periodic Witt vectors is
\begin{equation}
 \W_K^{(\mathfrak f)}(A)
 =\{x:\psi_{\mathfrak a}(x)=\psi_{\mathfrak b}(x)
          \text{ whenever }\mathfrak a\sim_{\mathfrak f}\mathfrak b\}
 \subseteq\W_K(A).
\label{eq:periodic-Witt}
\end{equation}
For flat \(A\), this says that the ghost coordinates factor through
\(\DR_{\mathfrak f}\).
\end{definition}

\noindent Put
\begin{equation}
 R_{\mathfrak f}:=\W_K^{(\mathfrak f)}(\Obar),\qquad
 R_K:=\varinjlim_{\mathfrak f}R_{\mathfrak f}.
\label{eq:Rf-periodic}
\end{equation}
For the ray class orders and their Witt-vector description, see
\cite[\S7.4, Proposition~9.8 and Corollary~9.10]{BdS18}.

\begin{theorem}
\label{thm:periodic-orders}
\noindent \(R_{\mathfrak f}\) is the finite flat maximal \(\Lambda_K\)-order
in \(E_{\mathfrak f}\), and
\begin{equation}
 K\tensor_{\Ok}R_{\mathfrak f}\cong E_{\mathfrak f},
 \qquad K\tensor_{\Ok}R_K\cong E_K.
\label{eq:global-generic}
\end{equation}
The identifications commute with conductor maps
\(u_{\mathfrak f,\mathfrak g}:R_{\mathfrak f}\into R_{\mathfrak g}\)
for \(\mathfrak f\mid\mathfrak g\).
If \(\mathfrak p\nmid\mathfrak f\), the localization
\((R_{\mathfrak f})_{\mathfrak p}\) is finite \(\et\) and is the
integral closure of \((\Ok)_{\mathfrak p}\) in \(E_{\mathfrak f}\).
As a filtered colimit of finite flat orders, \(R_K\) is flat over \(\Ok\).
\end{theorem}

\noindent In particular, we have the descriptions
\begin{equation}
 \DR_{\mathfrak f}\cong
 \Hom_{\Ok\text{-alg}}(R_{\mathfrak f},\Obar),\qquad
 \DR_K\cong\Hom_{\Ok\text{-alg}}(R_K,\Obar).
\label{eq:integral-Hom-points}
\end{equation}
Writing \(\beta_0\) for the first ghost coordinate,
\([\mathfrak a]\) corresponds to \(\beta_0\circ\psi_{\mathfrak a}\).
Thus multiplication by \(\mathfrak a\) on points is precomposition
by \(\psi_{\mathfrak a}\). Under \eqref{eq:global-generic} these are
the maps \(\sigma_{\mathfrak a}\) of Section~\ref{sec:DR}.
These are ordinary algebra homomorphisms, not
\(\Lambda_K\)-homomorphisms into \(\Obar\).
This description avoids choosing function coordinates for \(E_K\).

\begin{example}
\label{ex:cyclotomic-orders}
\noindent For \(K=\Q\), choose compatible primitive roots of unity
\(\zeta_{mn}^{\,m}=\zeta_n\).  Then
\begin{equation}
 R_n=\Z[T]/(T^n-1),\qquad R_{\Q}=\Z[\Q/\Z],
\label{eq:Q-global-order}
\end{equation}
with \(T\) corresponding to the periodic Teichm\"uller vector
\([\zeta_n]\), by \cite[Example~9.9]{BdS18}.
Write \(e(r)\) for the group-ring basis element at \(r\in\Q/\Z\).
All rational comparisons use this choice of roots of unity.
For example, \(R_2=\Z[T]/(T^2-1)\) has generic fibre
\(\Q\times\Q\), but it is not the product order:
the idempotent \((1+T)/2\) does not belong to \(R_2\), although
it is integral over \(\Z\).
This distinction is precisely what the next proposition detects.
\end{example}

\subsection{Why the partial inverses are not integral}
\label{sec:obstruction}

\begin{proposition}
\label{cor:rho-not-integral}
\noindent For every prime \(\mathfrak p\),
\(e_{\mathfrak p}\notin R_K\) and
\(\rho_{\mathfrak p}(R_K)\nsubseteq R_K\).
\end{proposition}

\begin{proof}
\noindent Suppose \(e_{\mathfrak p}\in B=R_{\mathfrak f}\).
Since \(\psi_{\mathfrak p}(e_{\mathfrak p})=1\) and
\(e_{\mathfrak p}^{\Nm(\mathfrak p)}=e_{\mathfrak p}\), the Frobenius
congruence gives \(1-e_{\mathfrak p}\in\mathfrak pB\).
The finite module \(C=(1-e_{\mathfrak p})B\) satisfies
\(C=\mathfrak pC\), so \(C_{\mathfrak p}=0\) by Nakayama.
Torsion freeness gives \(C=0\), contradicting the fact that
\(\mathfrak p\DR_K\) is proper.  Finally,
\(\rho_{\mathfrak p}(1)=e_{\mathfrak p}\).
\end{proof}

\begin{example}
\noindent Already over \(\Q\),
\begin{equation}
 \rho_n(e(r))=\frac1n\sum_{ns=r}e(s).
\label{eq:rho-Q-intro}
\end{equation}
Connes, Consani and Marcolli replace it by the integral additive map
\(\widetilde\rho_n=n\rho_n\) \cite[\S4]{CCM09}.
The denominator in this formula explains why the integral algebra
uses two asymmetric generators instead of a pair of integral
isometries.
\end{example}

\begin{remark}
\noindent The ordinary integral closure of \(R_K\) in \(E_K\) does contain
the idempotents and is stable under \(\rho\), but loses the inherited
Frobenius-lift structure: \(\sigma_{\mathfrak p}(e_{\mathfrak p})=1\)
is not congruent to \(e_{\mathfrak p}\) modulo \(\mathfrak p\), as
evaluation outside \(\mathfrak p\DR_K\) shows.
We instead work with \(R_K\) and Verschiebung.
\end{remark}

\subsection{Verschiebung}
\label{sec:V}

\noindent We now recall the Verschiebung family
\(V_{\mathfrak a,\xi}\) used in the introduction.  For an arbitrary \(\Ok\)-algebra \(A\),
Borger's construction gives natural \(\Ok\)-linear maps
\begin{equation}
 V_{\mathfrak a,\xi}:\W_K(A)\longrightarrow\W_K(A),
 \qquad \xi\in\mathfrak a,
\label{eq:V-a-xi}
\end{equation}
depending \(\Ok\)-linearly on \(\xi\).
They are constructed for prime powers in
\cite[\S4.2]{Borger11}. The general ideal version is described
in \cite[Remark~4.3]{Borger11}, using the iterated-Witt
decomposition of \cite[Corollary~5.4]{Borger11}.
For flat \(A\), they are uniquely characterized by
\begin{equation}
 (V_{\mathfrak a,\xi}x)_{\mathfrak b}
 =\begin{cases}
 \xi x_{\mathfrak b\mathfrak a^{-1}},&\mathfrak a\mid\mathfrak b,\\
 0,&\mathfrak a\nmid\mathfrak b.
 \end{cases}
\label{eq:V-ghost}
\end{equation}
For general \(A\), including residue fields, the operators are
defined by Borger's functorial construction, not by assuming
that the ghost map is injective.

\begin{remark}
\noindent Borger packages this bilinear family as a single map
\[
 \mathfrak a\tensor_{\Ok}\W_K(A)\longrightarrow\W_K(A),
 \qquad \xi\tensor x\longmapsto V_{\mathfrak a,\xi}(x).
\]
The tensor factor records the dependence on
\(\xi\in\mathfrak a\). It is not an additional Witt coordinate.
If \(\mathfrak a=(\alpha)\), every operator in the family is
an \(\Ok\)-multiple of \(V_{\mathfrak a,\alpha}\).
Thus a choice of generator recovers the familiar single operator.
We work with the family and express its interaction with
Frobenius by the projection formula below, without twisted-module
notation.
\end{remark}

\noindent We record the identities used below, without repeating
Borger's construction.  For \(x,z\in\W_K(A)\),
\(\xi\in\mathfrak a\),
\begin{equation}
 \psi_{\mathfrak a}V_{\mathfrak a,\xi}(x)=\xi x,\qquad
 V_{\mathfrak a,\xi}(\psi_{\mathfrak a}(z)x)
       =zV_{\mathfrak a,\xi}(x),
\label{eq:V-basic}
\end{equation}
\begin{equation}
 V_{\mathfrak a,\xi}V_{\mathfrak b,\eta}
       =V_{\mathfrak a\mathfrak b,\xi\eta}\quad(\eta\in\mathfrak b),
\qquad
 V_{\mathfrak a,\xi}(x)V_{\mathfrak a,\eta}(y)
       =V_{\mathfrak a,\xi\eta}(xy)
       \quad(\eta\in\mathfrak a).
\label{eq:V-composition}
\end{equation}
The last formula uses \(\mathfrak a^2\subseteq\mathfrak a\).
If \(\mathfrak d=\gcd(\mathfrak b,\mathfrak c)\), then
\begin{equation}
 \psi_{\mathfrak c}V_{\mathfrak b,\xi}
 =V_{\mathfrak b/\mathfrak d,\xi}\psi_{\mathfrak c/\mathfrak d},
 \qquad \xi\in\mathfrak b\subseteq\mathfrak b/\mathfrak d.
\label{eq:gcd-V}
\end{equation}
The projection and product formulas are
\cite[(4-2-4)-(4-2-6)]{Borger11}. The other identities follow
prime by prime from the same ghost principle and
\cite[Corollary~5.4]{Borger11}.
For \(K=\Q\) and \(\xi=b\), the ideal inclusion supplies the factor
\(\gcd(b,c)\) in the classical formula \(F_cV_b\).

\subsection{Restriction to periodic Witt vectors}

\begin{proposition}
\label{prop:V-periodic}
\noindent For \(\mathfrak a,\mathfrak f\in I_K\) and
\(\xi\in\mathfrak a\), one has
\begin{equation}
 V_{\mathfrak a,\xi}(R_{\mathfrak f})
       \subseteq R_{\mathfrak a\mathfrak f}.
\label{eq:V-periodic}
\end{equation}
These restrictions commute with conductor inclusions and hence
define an \(\Ok\)-linear map \(V_{\mathfrak a,\xi}:R_K\to R_K\).
Under \eqref{eq:global-generic}, it satisfies
\begin{equation}
 V_{\mathfrak a,\xi}(x)=\xi\rho_{\mathfrak a}(x)
 \qquad(x\in R_K).
\label{eq:V-is-xi-rho}
\end{equation}
\end{proposition}

\begin{proof}
\noindent Since \(\Obar\) is flat, its Witt ghost map is injective.
If \(\mathfrak b\sim_{\mathfrak a\mathfrak f}\mathfrak c\), the two
ideals have the same gcd with \(\mathfrak a\mathfrak f\).  Thus
\(\mathfrak a\) divides both or neither.  In the second case both
coordinates in \eqref{eq:V-ghost} vanish. In the first case,
\(\mathfrak b/\mathfrak a\sim_{\mathfrak f}\mathfrak c/\mathfrak a\)
by \eqref{eq:multiplication-ray-equivalence}, and periodicity of \(x\)
makes the coordinates equal.  This proves \eqref{eq:V-periodic}.
The same formula commutes with conductor inclusions.

\noindent The point \([\mathfrak b]\) in
\eqref{eq:integral-Hom-points} is the \(\mathfrak b\)-ghost coordinate.
Consequently \eqref{eq:V-ghost} and \eqref{eq:rho-functions} give
\eqref{eq:V-is-xi-rho} on every ideal point.  These points surject onto
each finite \(\DR_{\mathfrak f}\), hence separate elements of \(E_K\).
\end{proof}

\begin{remark}
\noindent The tensor product in Borger's formulation also has a
concrete interpretation here.  Since \(R_{\mathfrak f}\) is flat,
multiplication identifies
\[
 \mathfrak a\tensor_{\Ok}R_{\mathfrak f}
       \xrightarrow{\;\sim\;}\mathfrak aR_{\mathfrak f}
       \subset R_{\mathfrak f}\subset E_K
       \subset A_K^{\mathrm{arith}},
 \qquad \xi\tensor x\longmapsto\xi x.
\]
Under this identification, the map
\(\xi\tensor x\mapsto V_{\mathfrak a,\xi}(x)\)
is simply the restriction of \(\rho_{\mathfrak a}\) to
\(\mathfrak aR_{\mathfrak f}\).
The ordinary map \(\rho_{\mathfrak a}\) need not preserve
\(R_{\mathfrak f}\), but it sends this smaller ideal into
\(R_{\mathfrak a\mathfrak f}\).
We use the concrete ideals \(\mathfrak aR_K\) in the normal form
below.
\end{remark}


\section{Integral BC-algebra}
\label{sec:W-algebra}

\noindent We extend the integral BC-algebra for $\Q$ of
\cite[Definition~4.7 and Corollary~4.11]{CCM09} to general $K$.  We construct the algebra directly inside
\(A_K^{\mathrm{arith}}\), prove its normal form, and then extend
the Witt realization of \cite{CC14}.

\subsection{Construction and normal form of \(A_K^{\mathrm{int}}\)}

\noindent The elements \(U_{\mathfrak a}\) and
\(U_{\mathfrak a}^*\) are the generators of the rational crossed
product from Section~\ref{sec:DR}.

\begin{definition}
\label{def:W-integral-algebra}
\noindent The integral BC-algebra \(A_K^{\mathrm{int}}\) is the unital
\(\Ok\)-subalgebra of \(A_K^{\mathrm{arith}}\) generated by \(R_K\),
\(U_{\mathfrak a}^*\), and \(\xi U_{\mathfrak a}\), for
\(\mathfrak a\in I_K\) and \(\xi\in\mathfrak a\).
\end{definition}

\noindent The restriction \(\xi\in\mathfrak a\) is part of the
definition: the unscaled \(U_{\mathfrak a}\) need not be integral.
In particular, this is not the crossed product of \(R_K\) by
the maps \(\rho_{\mathfrak a}\).
The relations we use are the ordinary crossed-product relations
on these concrete elements:
\begin{equation}
\begin{gathered}
 U_{(1)}=U_{(1)}^*=1,\qquad
 U_{\mathfrak a}^*U_{\mathfrak b}^*=U_{\mathfrak a\mathfrak b}^*,\\
 (\xi U_{\mathfrak a})(\eta U_{\mathfrak b})
       =\xi\eta U_{\mathfrak a\mathfrak b},\qquad
 U_{\mathfrak a}^*(\xi U_{\mathfrak a})=\xi,
 \quad \xi\in\mathfrak a,\ \eta\in\mathfrak b,
\end{gathered}
\label{eq:relations-semigroup}
\end{equation}
and, for \(x\in R_K\),
\begin{equation}
\begin{gathered}
 U_{\mathfrak a}^*x=\psi_{\mathfrak a}(x)U_{\mathfrak a}^*,\qquad
 x(\xi U_{\mathfrak a})
       =(\xi U_{\mathfrak a})\psi_{\mathfrak a}(x),\\
 (\xi U_{\mathfrak a})xU_{\mathfrak a}^*=V_{\mathfrak a,\xi}(x).
\end{gathered}
\label{eq:relations-coeff}
\end{equation}
By Proposition~\ref{prop:V-periodic}, the product in the last equality
lies in \(R_K\).
For coprime ideals one also has
\begin{equation}
 (\xi U_{\mathfrak a})U_{\mathfrak b}^*
       =U_{\mathfrak b}^*(\xi U_{\mathfrak a})
       \qquad((\mathfrak a,\mathfrak b)=1).
\label{eq:relations-coprime}
\end{equation}
Define the \(\Ok\)-linear maps
\begin{equation}
 N_{\mathfrak a,\mathfrak b}:\mathfrak aR_K\longrightarrow
 A_K^{\mathrm{int}},\qquad
 z\longmapsto U_{\mathfrak a}zU_{\mathfrak b}^*.
\label{eq:normal-form-map}
\end{equation}
They take values in the integral algebra because
\(z=\sum_i\xi_i x_i\) gives
\(U_{\mathfrak a}zU_{\mathfrak b}^*
=\sum_i(\xi_iU_{\mathfrak a})x_iU_{\mathfrak b}^*\).

\begin{theorem}
\label{thm:witt-integral-algebra}
\noindent The algebra has the internal direct-sum decomposition
\begin{equation}
 A_K^{\mathrm{int}}
 =\bigoplus_{\substack{\mathfrak a,\mathfrak b\in I_K\\
                      (\mathfrak a,\mathfrak b)=1}}
       U_{\mathfrak a}(\mathfrak aR_K)U_{\mathfrak b}^*.
\label{eq:module-normal-form}
\end{equation}
Each map \(N_{\mathfrak a,\mathfrak b}\) is injective.
In particular, \(A_K^{\mathrm{int}}\) is flat over \(\Ok\), and
\begin{equation}
 K\tensor_{\Ok}A_K^{\mathrm{int}}
       \xrightarrow{\;\sim\;}E_K\rtimes_\rho I_K.
\label{eq:W-generic-fibre}
\end{equation}
\end{theorem}

\begin{proof}
\noindent We first reduce words in the concrete generators.
If \(\mathfrak d=\gcd(\mathfrak c,\mathfrak b)\), put
\(\mathfrak c=\mathfrak d\mathfrak c_0\) and
\(\mathfrak b=\mathfrak d\mathfrak b_0\).
Cancellation and coprimality in the rational crossed product give
\begin{equation}
 U_{\mathfrak c}^*(\xi U_{\mathfrak b})
       =(\xi U_{\mathfrak b_0})U_{\mathfrak c_0}^*.
\label{eq:mixed-reduction}
\end{equation}
This is still an integral expression, since
\(\xi\in\mathfrak b\subseteq\mathfrak b_0\).
Moving coefficients by \eqref{eq:relations-coeff} and applying
\eqref{eq:mixed-reduction} moves all adjoint generators to the
right.  Hence every word is a sum of
\((\xi U_{\mathfrak a})xU_{\mathfrak b}^*\).

\noindent To make the indices coprime, write
\(\mathfrak a=\mathfrak d\mathfrak a'\),
\(\mathfrak b=\mathfrak d\mathfrak b'\), where
\(\mathfrak d=\gcd(\mathfrak a,\mathfrak b)\).
Because \(\xi\in\mathfrak d\mathfrak a'\), write
\(\xi=\sum_i\delta_i\alpha_i\), with
\(\delta_i\in\mathfrak d\), \(\alpha_i\in\mathfrak a'\).
Then
\begin{equation}
 (\xi U_{\mathfrak a})xU_{\mathfrak b}^*
 =\sum_i(\alpha_iU_{\mathfrak a'})
       V_{\mathfrak d,\delta_i}(x)U_{\mathfrak b'}^*.
\label{eq:cancel-normal}
\end{equation}
Verschiebung preserves \(R_K\), so these expressions prove spanning.

\noindent Different coprime pairs have different \(J_K\)-degrees.
Within one degree, the coefficient is recovered by
\[
 U_{\mathfrak a}^*
       (U_{\mathfrak a}zU_{\mathfrak b}^*)U_{\mathfrak b}=z.
\]
Thus each \(N_{\mathfrak a,\mathfrak b}\) is injective and the
sum is direct.  Each \(\mathfrak aR_K\) is torsion free over the
Dedekind domain \(\Ok\), hence flat, and so is their direct sum.
Finally, \(K\tensor_{\Ok}R_K=E_K\), and inverting any
\(0\ne\xi\in\mathfrak a\) recovers
\(U_{\mathfrak a}=\xi^{-1}(\xi U_{\mathfrak a})\).
Consequently the generic fibre contains all crossed-product
generators, proving \eqref{eq:W-generic-fibre}.
\end{proof}


\noindent The normal form specifies which coefficients belong
to the integral model in every ideal degree.

\begin{corollary}
\label{cor:concrete-integral}
\noindent Write \(\gamma=\mathfrak a\mathfrak b^{-1}\) in coprime
form and put
\(A_{K,\gamma}^{\mathrm{int}}
=A_K^{\mathrm{int}}\cap A_{K,\gamma}^{\mathrm{arith}}\).
Then
\[
 A_{K,\gamma}^{\mathrm{arith}}=U_{\mathfrak a}E_KU_{\mathfrak b}^*,
 \qquad
 A_{K,\gamma}^{\mathrm{int}}
       =U_{\mathfrak a}(\mathfrak aR_K)U_{\mathfrak b}^*.
\]
For \(z\in A_{K,\gamma}^{\mathrm{arith}}\), the membership criterion is
\[
 z\in A_K^{\mathrm{int}}
 \quad\Longleftrightarrow\quad
 U_{\mathfrak a}^*zU_{\mathfrak b}\in\mathfrak aR_K.
\]
In particular, \(A_K^{\mathrm{int}}\cap E_K=R_K\).
\end{corollary}

\begin{proof}
\noindent The rational normal form writes
\(z=U_{\mathfrak a}xU_{\mathfrak b}^*\), and compression recovers
\(x\).  Theorem~\ref{thm:witt-integral-algebra} gives the criterion.
Degree \(1\) corresponds to \(\mathfrak a=\mathfrak b=(1)\).
\end{proof}

\begin{remark}
\noindent Under the identification of the degree
\(\mathfrak a\mathfrak b^{-1}\) component with \(\mathfrak aR_K\),
left and right multiplication by \(r,s\in R_K\) act as
\[
 z\longmapsto\psi_{\mathfrak a}(r)z\psi_{\mathfrak b}(s).
\]
The grading is not strong: products of homogeneous components satisfy
\[
 A_{K,\mathfrak a^{-1}}^{\mathrm{int}}
 A_{K,\mathfrak a}^{\mathrm{int}}=\mathfrak aR_K,\qquad
 A_{K,\mathfrak a}^{\mathrm{int}}
 A_{K,\mathfrak a^{-1}}^{\mathrm{int}}
       =\rho_{\mathfrak a}(\mathfrak aR_K).
\]
These equalities follow from \eqref{eq:relations-coeff}. Their right
sides need not equal \(R_K\).  In particular, the concrete realization
does not make \(A_K^{\mathrm{int}}\) a crossed product by
endomorphisms \(\rho_{\mathfrak a}\) of \(R_K\).
\end{remark}

\begin{corollary}
\label{thm:Q-specialization}
\noindent For \(K=\Q\), the compatible roots chosen in
Section~\ref{sec:orders} identify \(A_{\Q}^{\mathrm{int}}\)
with the integral BC-algebra of \cite[Definition~4.7]{CCM09}, by
\begin{equation}
 R_{\Q}=\Z[\Q/\Z],\qquad
 U_{(n)}^*=\mu_n^*,\qquad
 nU_{(n)}=\widetilde\mu_n,\qquad
 V_{(n),n}=n\rho_n.
\label{eq:Q-identification}
\end{equation}
The last operator is denoted \(\widetilde\rho_n\) in
\cite[\S4.3]{CCM09}.
\end{corollary}

\begin{proof}
\noindent Every \(\xi\in(n)\) is an integer multiple of \(n\).
Equations \eqref{eq:Q-global-order} and \eqref{eq:V-ghost} identify
the commutative ring and Verschiebung operators.
The generators and normal form are then those of \cite[\S4]{CCM09}.
\end{proof}

\begin{remark}
\noindent The usual BC-involution need not preserve the integral
algebra: over \(\Q\), it sends \(U_{(n)}^*\) to \(U_{(n)}\),
which does not belong to \(A_{\Q}^{\mathrm{int}}\) for \(n>1\).
Our integral construction is an algebraic one.
\end{remark}

\subsection{Witt-vector representations}
\label{subsec:Witt-representation}

\noindent The representation of \cite[Theorems~4.5 and 6.4]{CC14}
realizes the integral operators as actual Frobenius and
Verschiebung on Witt vectors in characteristic \(p\).  This
provides an arithmetic realization of the algebra, distinct from
the scalar extension in Part~\ref{part:padic}.
Fix a rational prime \(p\) and an embedding
\(\jmath:\Kbar\into\C_p\).

Put
\[
 \mathfrak m_{\C_p}:=\{z\in\Ocal_{\C_p}:|z|_p<1\},
 \qquad
 \overline k:=\Ocal_{\C_p}/\mathfrak m_{\C_p}.
\]
This residue field is an algebraic closure of \(\F_p\).
Algebraic integers map into \(\Ocal_{\C_p}\), so \(\jmath\)
induces
\[
 \Obar\longrightarrow\Ocal_{\C_p}
          \longrightarrow\overline k.
\]
The kernel of the composite intersects \(\Ok\) in a prime
\(\mathfrak p\mid p\), determined by \(\jmath\).
Give \(\overline k\) the resulting \(\Ok\)-algebra structure.  The induced map on Witt vectors is
\[
 \red:\W_K(\Obar)\longrightarrow \W_K(\overline k).
\]

\begin{proposition}
\label{thm:witt-vector-representation}
\noindent There is a unique unital \(\Ok\)-algebra homomorphism
\[
 \Pi:A_K^{\mathrm{int}}\longrightarrow
 \operatorname{End}_{\Ok}\bigl(\W_K(\overline k)\bigr)
\]
with the following action on generators. For \(x\in R_K\),
\(\mathfrak a\in I_K\), \(\xi\in\mathfrak a\), and
\(w\in\W_K(\overline k)\),
\begin{equation}
 \Pi(x)w=\red(x)w,\qquad
 \Pi(U_{\mathfrak a}^*)w=\psi_{\mathfrak a}(w),\qquad
 \Pi(\xi U_{\mathfrak a})w=V_{\mathfrak a,\xi}(w).
\label{eq:Pi-formulas}
\end{equation}
\end{proposition}

\begin{proof}
\noindent We define the action on the unique normal form.
For a coprime pair \((\mathfrak a,\mathfrak b)\) and
\(z=\sum_i\xi_i x_i\in\mathfrak aR_K\), with
\(\xi_i\in\mathfrak a\), \(x_i\in R_K\), put
\[
 \Pi(U_{\mathfrak a}zU_{\mathfrak b}^*)(w)
 :=\sum_i V_{\mathfrak a,\xi_i}
       \bigl(\red(x_i)\psi_{\mathfrak b}(w)\bigr).
\]
Each summand is additive in \(\xi_i\) and \(x_i\). Moving a scalar
\(c\in\Ok\) from \(\xi_i\) to \(x_i\) does not change it:
\(V_{\mathfrak a,c\xi}(y)=cV_{\mathfrak a,\xi}(y)
=V_{\mathfrak a,\xi}(cy)\), and \(\red(cx)=c\red(x)\).
Thus the expression factors through
\(\mathfrak a\tensor_{\Ok}R_K\simeq\mathfrak aR_K\)
and depends only on \(z\).
The direct-sum normal form therefore defines a linear map
\(\Pi\) on the whole algebra.

\noindent To check multiplicativity, apply the word reductions
in the proof of Theorem~\ref{thm:witt-integral-algebra}.
Frobenius and Verschiebung satisfy the semigroup and covariance
rules by \eqref{eq:V-basic}--\eqref{eq:V-composition}.
The mixed reduction \eqref{eq:mixed-reduction} is exactly
\eqref{eq:gcd-V}.  For the cancellation
\eqref{eq:cancel-normal}, the projection formula gives
\[
 V_{\mathfrak d,\delta}
       \bigl(\red(x)\psi_{\mathfrak d}(w)\bigr)
       =V_{\mathfrak d,\delta}(\red(x))\,w
       =\red(V_{\mathfrak d,\delta}(x))\,w,
\]
where the last equality is naturality of Verschiebung.
Thus each reduction of a product agrees with composition
of the corresponding operators.  This proves multiplicativity
and the formulas \eqref{eq:Pi-formulas}.
\end{proof}

\begin{example}
\noindent Over \(\Q\), this is the representation of
\cite[Theorems~4.5 and 6.4]{CC14}:
\(\mu_n^*\) acts by \(F_n\) and \(\widetilde\mu_n\) by \(V_n\).
Indeed, on coefficients it is the composite
\[
 \Z[\Q/\Z]\longrightarrow\Z[\mu^{(p)}]
 \cong\W_0(\overline k)\lhook\joinrel\longrightarrow
 \W(\overline k).
\]
Here \(\mu^{(p)}\) is the group of roots of unity of order prime to
\(p\). The first map removes the \(p\)-primary part, and the last
map is the Almkvist embedding of the rational Witt ring.
For example, the integral generator acts by
\(V_p=pF_p^{-1}\), since \(F_p\) is invertible over the perfect field
\(\overline k\).
\end{example}

\begin{remark}
\noindent The inclusion \(A_K^{\mathrm{int}}\subseteq A_K^{\mathrm{arith}}\)
does not imply that \(\Pi\) is faithful.
Over \(\Q\),
\(e(1/p)-1\ne0\) in \(R_{\Q}\), but its Witt reduction is zero
because a \(p\)-th root of unity reduces to \(1\).
The \(p\)-adic KMS forms below are constructed on the
scalar-extended algebra, not by requiring descent through this
representation.
\end{remark}

\noindent Over \(\Q\), Connes and Consani identify the group ring on the
roots of unity of order prime to \(p\) with \(W_0(\overline k)\)
inside \(W(\overline k)\) \cite[Theorem~2.1 and Proposition~4.4]{CC14}.
The following proposition gives the corresponding injectivity statement
for conductors prime to \(\mathfrak p\) over \(K\).

\begin{proposition}
\label{prop:witt-reduction-injective}
\noindent Put
\[
 R_{K,\mathfrak p}:=
       \varinjlim_{\mathfrak p\nmid\mathfrak f}R_{\mathfrak f}.
\]
Then \(\red:R_{K,\mathfrak p}\to \W_K(\overline k)\) is injective.
Its image is stable under all Frobenius operators and under
\(V_{\mathfrak a,\xi}\) for \(\mathfrak p\nmid\mathfrak a\)
and \(\xi\in\mathfrak a\).
In particular, \(\Pi\) is faithful on \(R_{K,\mathfrak p}\).
\end{proposition}

\begin{proof}
\noindent Let \(A=(\Ok)_{\mathfrak p}\), choose a uniformizer
\(\pi\) of \(A\), and write \(\W_{\mathfrak p}\) for its
single-prime Witt functor.
Localization and decomposition of Witt vectors
\cite[\S1.20 and Corollary~5.4]{Borger11} give
\[
 \W_K(\overline k)\simeq
 \prod_{\substack{\mathfrak b\in I_K\\\mathfrak p\nmid\mathfrak b}}
       \W_{\mathfrak p}(\overline k).
\]
The projection to the factor indexed by \(\mathfrak b\) applies
\(\psi_{\mathfrak b}\) and then restricts to the
\(\mathfrak p\)-typical coordinates.
In each factor Frobenius is bijective because \(\overline k\)
is perfect, Verschiebung \(V_\pi\) is an injective shift, and
\(FV_\pi=\pi\)
\cite[\S\S3.5,3.7 and (4-2-5)]{Borger11}.
Thus \(\W_K(\overline k)\) is an \(A\)-algebra with no \(\pi\)-torsion.
No unramifiedness of \(K_{\mathfrak p}/\Q_p\) is needed.

\noindent Fix \(\mathfrak p\nmid\mathfrak f\).
The algebra \(B=(R_{\mathfrak f})_{\mathfrak p}\) is finite
\'etale over \(A\)
\cite[proof of Proposition~7.3(1)]{BdS18}, and Witt reduction
extends to \(B\to \W_K(\overline k)\).
If \(x\) lies in its kernel, naturality shows that every
evaluation \(\beta_0\psi_{\mathfrak b}(x)\), for
\(\mathfrak p\nmid\mathfrak b\), vanishes after reduction.
These ideals represent all points of \(\DR_{\mathfrak f}\):
in each ray-class stratum an ideal representative may be chosen
away from \(\mathfrak p\).
After base change to \(\Ocal_{\C_p}\), the finite \'etale algebra
\(B\) splits as a product of copies of \(\Ocal_{\C_p}\), since
\(\Ocal_{\C_p}\) is strictly henselian. Hence reducing the generic
evaluations gives every geometric point of \(B/\pi B\).
Consequently \(x\in\pi B\).
Write \(x=\pi y\).  Since \(\W_K(\overline k)\) has no \(\pi\)-torsion,
\(\red(y)=0\). Iteration gives
\(x\in\bigcap_{n\ge0}\pi^nB=0\).
Localization is injective on the torsion-free module
\(R_{\mathfrak f}\), so passing to the colimit proves injectivity.

\noindent Stability follows from Frobenius periodicity,
Proposition~\ref{prop:V-periodic}, and naturality.
Finally, multiplication by \(\red(x)\) on \(\W_K(\overline k)\) recovers
\(\red(x)\) at \(1\in \W_K(\overline k)\), proving faithfulness on
\(R_{K,\mathfrak p}\).
\end{proof}


\part{\texorpdfstring{\(p\)}{p}-adic BC-systems}
\label{part:padic}

\section{The \texorpdfstring{\(p\)}{p}-adic BC-system}
\label{sec:padic}

\noindent We follow the algebraic \(p\)-adic viewpoint of
Connes and Consani \cite[\S\S6-7]{CC14}.
We first distinguish the scalar-extended algebra from its integral
model, then define the automorphisms and the \(p\)-adic KMS condition.
The existence of the relevant linear forms is proved in
Section~\ref{sec:DR-padic}, using the results of Deligne-Ribet \cite{DR80}.

\subsection{The $p$-adic BC-algebra}

\noindent Fix a prime number \(p\) and an embedding \(\jmath:\Kbar\into\C_p\).
We also write \(\jmath\) for its restriction to \(K\).
As in the rational construction of Connes and Consani
\cite[Definition~6.7]{CC14}, we restrict the ideal operators to
indices prime to \(p\). Unlike their construction, we use the
whole commutative algebra \(R_K\), without taking a quotient.
Recall that
\begin{equation}
 I_K^{(p)}:=\{\mathfrak a\in I_K:(\mathfrak a,p\Ok)=1\},
 \qquad J_K^{(p)}:=\langle I_K^{(p)}\rangle\subseteq J_K,
\label{eq:prime-to-p-monoid}
\end{equation}
and let \(A_K^{\mathrm{int},(p)}\subseteq A_K^{\mathrm{int}}\) be generated by \(R_K\) and
\(U_{\mathfrak a}^*\) and \(\xi U_{\mathfrak a}\), for
\(\mathfrak a\in I_K^{(p)}\) and \(\xi\in\mathfrak a\).  Define
\begin{equation}
 A_{K,p}:=\C_p\tensor_{\Ok,\jmath}A_K^{\mathrm{int},(p)},
 \qquad E_{K,p}:=\C_p\tensor_{\Ok,\jmath}R_K
       \cong\C_p\tensor_{K,\jmath}E_K.
\label{eq:padic-coefficient-algebra}
\end{equation}
If \(\Ocal_{\C_p}\) is the valuation
ring, the original integral model gives the subalgebra
\(\Ocal_{\C_p}\tensor_{\Ok,\jmath}A_K^{\mathrm{int},(p)}
\subseteq A_{K,p}\), by flatness and the normal-form theorem.
This scalar extension is distinct from the residue-Witt
representation of Section~\ref{subsec:Witt-representation}.

\begin{proposition}
\label{prop:padic-generic}
\noindent There are canonical graded isomorphisms
\[
 A_{K,p}\cong
 \C_p\tensor_{K,\jmath}
          (E_K\rtimes_\rho I_K^{(p)})
 \cong E_{K,p}\rtimes_\rho I_K^{(p)}.
\]
The generators \(U_{\mathfrak a}\) and \(U_{\mathfrak a}^*\)
are now available without an integrality restriction:
for \(0\ne\xi\in\mathfrak a\), the image of
\(\xi U_{\mathfrak a}\) can be divided by \(\jmath(\xi)\).
\end{proposition}

\begin{proof}
\noindent Apply the normal-form theorem to the coprime pairs in
\(I_K^{(p)}\), then extend scalars.  Each ideal coefficient
\(\mathfrak a\tensor_{\Ok}\C_p\) becomes \(\C_p\), and the
normalized generators satisfy \eqref{eq:crossed-relations}.
\end{proof}

\begin{remark}
\label{rem:integral-versus-generic}
\noindent Thus the algebra over \(\C_p\) can be defined directly
from the rational endomotive: it does not depend on first
constructing a particular integral order. The periodic Witt vector model specifies the integral subalgebra, the
Frobenius-lift structure, the integral
Verschiebung and residue-Witt operators.
 The Deligne-Ribet KMS forms below live on the generic
 fibre and are compatible with its integral operators.
\end{remark}

\subsection{\texorpdfstring{\(p\)}{p}-adic time evolution}
\label{subsec:padic-dynamics}

\noindent We adapt the automorphisms of Connes and Consani
\cite[Lemmas~7.1--7.2 and Proposition~7.3]{CC14} to the ideal norms of \(K\).
Put \(q_p=p\) for odd \(p\), and \(q_2=4\), and normalize
\(v_p\) by \(v_p(p)=1\).
Write \(u=\omega(u)\langle u\rangle\) for \(u\in\Z_p^\times\),
where \(\langle u\rangle\in1+q_p\Z_p\).
Here \(\omega\) is the Teichm\"uller character for odd \(p\) and
the sign modulo \(4\) for \(p=2\).
On these units the Iwasawa logarithm is given by
\(\log_p u=\log_p\langle u\rangle\), with
\(\log_p(1+t)=\sum_{n\ge1}(-1)^{n+1}t^n/n\).
It kills the finite factor \(\omega(u)\). Define
\begin{equation}
 H_{K,p}:=\overline{\Nm(J_K^{(p)})}\subseteq\Z_p^\times,\qquad
 \log_p\langle H_{K,p}\rangle=p^{m_{K,p}}\Z_p.
\label{eq:norm-subgroup}
\end{equation}
The closure is taken in \(\Z_p^\times\).
The logarithmic image is a closed additive subgroup of \(\Z_p\).
It is nonzero and \(m_{K,p}\) is finite: if \(d=[K:\Q]\), norms of
rational principal ideals give
\((\Z_p^\times)^d\subseteq H_{K,p}\), and hence
\[
 v_p(q_p)\le m_{K,p}\le v_p(q_p)+v_p(d).
\]
Set
\begin{equation}
 D_{K,p}:=\{\beta\in\C_p:
       |\beta|_p<p^{m_{K,p}-1/(p-1)}\}.
\label{eq:Iwasawa-disk}
\end{equation}
For \(u\in H_{K,p}\) and \(\beta\in D_{K,p}\), put
\begin{equation}
 u^{(\beta)}
 :=\omega(u)\exp\bigl(\beta\log_p\langle u\rangle\bigr)
 =u\exp\bigl((\beta-1)\log_p u\bigr).
\label{eq:Iwasawa-character}
\end{equation}
The logarithms have absolute value at most \(p^{-m_{K,p}}\).
Thus the exponential converges uniformly in \(u\) on closed
subdiscs of \(D_{K,p}\), is multiplicative in \(u\), and is
analytic in \(\beta\).  The radius is greater than \(1\), so
\(D_{K,p}\) contains \(\Z_p\) and is stable under \(\beta\mapsto1-\beta\).
On this enlarged disc, \eqref{eq:Iwasawa-character} is asserted on
\(H_{K,p}\), not on all of \(\Z_p^\times\).

\begin{definition}
\label{def:padic-system}
\noindent The \emph{\(p\)-adic BC-system} is \(A_{K,p}\) with the
family of \(\C_p\)-algebra automorphisms \(\sigma^{(\beta)}\),
\(\beta\in D_{K,p}\). For each homogeneous \(A_\gamma\) of degree
\(\gamma\in J_K^{(p)}\), set
\[
 \sigma^{(\beta)}(A_\gamma)
       :=\bigl((\Nm(\gamma))^{(\beta)}\bigr)^{-1}A_\gamma.
\]
\end{definition}

\noindent Multiplicativity makes each map an algebra automorphism.
Equivalently,
\begin{equation}
 \sigma^{(\beta)}(U_{\mathfrak a}xU_{\mathfrak b}^*)
 =\left(\frac{\Nm(\mathfrak b)}{\Nm(\mathfrak a)}\right)^{(\beta)}
       U_{\mathfrak a}xU_{\mathfrak b}^*
 \qquad(x\in E_{K,p}).
\label{eq:p-adic-dynamics}
\end{equation}
These automorphisms play the role of imaginary-time evolution
in the \(p\)-adic KMS condition. Since the finite character is fixed,
\[
 \sigma^{(\beta)}\sigma^{(\gamma)}
       =\sigma^{(\beta+\gamma)}\sigma^{(0)},\qquad
 \sigma^{(0)}(A_\delta)=\omega(\Nm(\delta))^{-1}A_\delta.
\]
Thus it is \(\beta\mapsto\sigma^{(\beta)}(\sigma^{(0)})^{-1}\)
that is an additive flow on \(D_{K,p}\). We keep the finite
factor in the KMS condition to match Deligne-Ribet interpolation.

\begin{remark}
\label{rem:parameter-disc}
\noindent For \(K=\Q\), \(H_{\Q,p}=\Z_p^\times\) and
\(m_{\Q,p}=v_p(q_p)\), so \(D_{\Q,p}\) is exactly the
Connes-Consani extended \(s\)-disc.
The bound above also gives \(D_{K,p}=D_{\Q,p}\) whenever
\(p\nmid[K:\Q]\).
For every \(K\), if \(\varphi(q_p)\mid k\), then
\(u^{(1-k)}=u^{1-k}\) on \(H_{K,p}\).
At \(p=2\), the finite factor is the sign modulo \(4\) and
cannot in general be dropped.
\end{remark}

\begin{remark}
\label{rem:parameter-cover}
\noindent Connes and Consani go beyond the exponential disc by
using the multiplicative group
\(\{t\in\C_p:|t-1|_p<1\}\), whose logarithm covers \(\C_p\)
with kernel the \(p\)-power roots of unity
\cite[\S7.5]{CC14}.  This distinguishes characters invisible
to the logarithmic parameter \(\beta\).
Our enlargement is still an exponential disc, not a replacement
for that covering.  It suffices for the interpolation and KMS
construction proved here.  Extending the full construction over
the covering would require keeping those finite-order choices
and reexamining the smoothing denominators.
\end{remark}

\subsection{The \texorpdfstring{\(p\)}{p}-adic KMS condition}

\noindent On analytic elements of a complex BC-system, the usual
KMS condition can be written
\(\varphi(X\sigma_{i\beta}(Y))=\varphi(YX)\).
Since \(\sigma_t(U_{\mathfrak a})=\Nm(\mathfrak a)^{it}
U_{\mathfrak a}\), imaginary time scales \(U_{\mathfrak a}\) by
\(\Nm(\mathfrak a)^{-\beta}\).
Our automorphisms replace this scalar by
\(\bigl((\Nm(\mathfrak a))^{(\beta)}\bigr)^{-1}\), using
\eqref{eq:Iwasawa-character}.
The resulting identity is an algebraic twisted-trace condition:
it makes sense for any algebra with a specified automorphism.  Following \cite[\S7.4]{CC14}, we use it without
imposing the positivity or strip-holomorphy axioms of complex
\(C^*\)-theory.

\begin{definition}
\label{def:padic-KMS}
\noindent Let \(\beta\in D_{K,p}\).
A \(\C_p\)-linear form \(\varphi:A_{K,p}\to\C_p\) satisfies the
\emph{\(p\)-adic KMS condition at \(\beta\)} if
\begin{equation}
 \varphi\bigl(X\sigma^{(\beta)}(Y)\bigr)=\varphi(YX)
 \qquad(X,Y\in A_{K,p}).
\label{eq:padic-KMS}
\end{equation}
It is a \emph{\(p\)-adic KMS state} if also \(\varphi(1)=1\).
No positivity condition is imposed.
\end{definition}

\noindent For a meromorphic family \(\varphi_\beta\), its
\emph{partition function} is \(Z(\beta):=\varphi_\beta(1)\).
This refers to the chosen unnormalized family. The KMS condition
alone does not fix its scalar normalization.
Where the family is regular and \(Z(\beta)\ne0\),
\(Z(\beta)^{-1}\varphi_\beta\) is a \(p\)-adic KMS state.

\noindent We consider unnormalized forms as well because their values will
be $p$-adic \(L\)-values, whereas normalization gives ratios.  At a zero of \(Z\), extension of the normalized family requires
a separate argument.  In the real-quadratic case, the pole at
\(\beta=1\) is handled by Corollary~\ref{cor:quadratic-weight-one}.

\noindent The normal form of Theorem~\ref{thm:witt-integral-algebra},
after scalar extension as in Proposition~\ref{prop:padic-generic},
is the direct sum
\begin{equation}
 A_{K,p}
 =\bigoplus_{\substack{\mathfrak a,\mathfrak b\in I_K^{(p)}\\
                      (\mathfrak a,\mathfrak b)=1}}
       U_{\mathfrak a}E_{K,p}U_{\mathfrak b}^*.
\label{eq:padic-normal-form}
\end{equation}
Thus every element has a unique finite expression
\(\sum U_{\mathfrak a}x_{\mathfrak a,\mathfrak b}U_{\mathfrak b}^*\),
with \(x_{\mathfrak a,\mathfrak b}\in E_{K,p}\) and coprime indices.
Each summand has fractional-ideal degree
\(\mathfrak a\mathfrak b^{-1}\). Degree \(1\) occurs only for
\(\mathfrak a=\mathfrak b=(1)\).
We therefore have the linear projection
\begin{equation}
 \mathbb E_0:A_{K,p}\longrightarrow E_{K,p},\qquad
 \mathbb E_0\left(\sum
 U_{\mathfrak a}x_{\mathfrak a,\mathfrak b}U_{\mathfrak b}^*\right)
       =x_{(1),(1)}.
\label{eq:projection-abelian-part}
\end{equation}
The KMS condition for forms
\(\ell\circ\mathbb E_0\) reduces to the following identity for
the scalar-extended Verschiebung operators
\(V_{\mathfrak a,\xi}\) of \eqref{eq:V-is-xi-rho}.

\begin{proposition}
\label{prop:KMS-criterion}
\noindent Fix \(\beta\in D_{K,p}\).
For a \(\C_p\)-linear form \(\ell:E_{K,p}\to\C_p\), the form
\(\ell\circ\mathbb E_0\) satisfies \eqref{eq:padic-KMS}
if and only if
\begin{equation}
 \ell\bigl(V_{\mathfrak a,\xi}(x)\bigr)
 =\jmath(\xi)
   \bigl((\Nm(\mathfrak a))^{(\beta)}\bigr)^{-1}\ell(x)
\label{eq:distribution-required}
\end{equation}
for \(\mathfrak a\in I_K^{(p)}\), \(\xi\in\mathfrak a\), and
\(x\in E_{K,p}\).
\end{proposition}

\begin{proof}
\noindent For necessity, take \(X=(\xi U_{\mathfrak a})x\) and
\(Y=U_{\mathfrak a}^*\).  The crossed-product relations give
\((\Nm(\mathfrak a))^{(\beta)}
\ell(V_{\mathfrak a,\xi}(x))=\jmath(\xi)\ell(x)\).
For sufficiency, choose any \(0\ne\xi\in\mathfrak a\). Dividing
\eqref{eq:distribution-required} by \(\jmath(\xi)\) gives
\(\ell(\rho_{\mathfrak a}x)=((\Nm(\mathfrak a))^{(\beta)})^{-1}\ell(x)\).
It suffices to check homogeneous monomials. Both sides of
\eqref{eq:padic-KMS} vanish unless their degrees are inverse.
In the remaining case write
\[
 X=U_{\mathfrak a}xU_{\mathfrak b}^*,\qquad
 Y=U_{\mathfrak b}yU_{\mathfrak a}^*,\qquad
 (\mathfrak a,\mathfrak b)=1.
\]
Now \(XY=\rho_{\mathfrak a}(xy)\) and \(YX=\rho_{\mathfrak b}(xy)\).
Using \eqref{eq:p-adic-dynamics} for \(Y\), both sides of the KMS
condition equal \(((\Nm(\mathfrak b))^{(\beta)})^{-1}\ell(xy)\).
\end{proof}

\begin{remark}
\label{rem:diagonal-KMS}
\noindent In this paper, we only construct forms of the form
\(\ell \circ \mathbb E_0\). We do not classify all \(p\)-adic KMS forms.
\end{remark}

\subsection{Comparison with the rational case}

\noindent Over \(\Q\), our algebra maps onto the one used by
Connes and Consani.  This comparison recovers their KMS forms
without identifying two different commutative algebras.

\begin{example}
\noindent Recall that
\(\widetilde\rho_n:=n\rho_n=V_{(n),n}\). On the group-ring basis,
\[
 \widetilde\rho_n(e(r))=\sum_{ns=r}e(s).
\]
For \(p\nmid n\), the criterion of
Proposition~\ref{prop:KMS-criterion} becomes
\[
 \ell(\widetilde\rho_n(x))
 =\frac{n}{n^{(\beta)}}\ell(x)
 =\langle n\rangle^{1-\beta}\ell(x).
\]
Connes and Consani verify this formula on their quotient algebra
in \cite[Lemma~7.10]{CC14}.

\noindent Our commutative algebra is \(\C_p[\Q/\Z]\), whereas
\cite[Proposition~6.6 and Definition~6.7]{CC14} uses
\(\C_p[(\Q/\Z)^{(p)}]\), where \((\Q/\Z)^{(p)}\) consists of
elements of order prime to \(p\).  Witt reduction kills
\(p\)-power roots of unity, giving the quotient
\[
 \pi:\C_p[\Q/\Z]\onto\C_p[(\Q/\Z)^{(p)}],
 \qquad e(r)\longmapsto e(r_{p'}),
\]
where \(r_{p'}\) is the prime-to-\(p\) part of \(r\).
For instance, \(e(1/p)-1\) is a nonzero element of its kernel.
For \(p\nmid n\), \(\pi\) commutes with Frobenius and
\(\widetilde\rho_n\), hence extends to the crossed products
and intertwines their automorphisms.
The Connes-Consani KMS forms therefore pull back to
\(A_{\Q,p}\), with the same Kubota-Leopoldt partition function
\cite[Theorem~7.11]{CC14}.
\end{example}

\noindent Thus the rational comparison is a recovery by quotient,
not an equality of algebras.  For totally real fields, the next
section constructs forms on \(\C_p\tensor_K E_K\) directly from
Deligne-Ribet measures, without imposing a residue-Witt quotient.

\subsection{Arithmetic symmetries}

\noindent Identify \(G_K^{\ab}=\Gal(K^{\ab}/K)\) with
\(\DR_K^\times\).  For \(g\in G_K^{\ab}\), multiplication
\(x\mapsto gx\) is a Galois-equivariant bijection of
\(\DR_K=\Hom_{K\text{-alg}}(E_K,\Kbar)\).
The corresponding \(K\)-algebra automorphism
\(\theta_g:E_K\to E_K\) is characterized by
\[
 x(\theta_gF)=(gx)(F)
 \qquad(F\in E_K,\ x\in\DR_K).
\]
In evaluation notation this says \((\theta_gF)(x)=F(gx)\),
where \(gx\) is the monoid product.
Extend this action \(\C_p\)-linearly to \(E_{K,p}\):
\(\theta_g(\lambda\tensor F)=\lambda\tensor\theta_g(F)\).
It does not act on the scalar \(\lambda\in\C_p\).

\begin{proposition}
\label{prop:padic-symmetry}
\noindent For \(g\in G_K^{\ab}\), the action \(\theta_g\)
preserves \(R_K\), commutes with Frobenius and Verschiebung, and naturally 
extends to \(A_K^{\mathrm{int}}\) and \(A_{K,p}\).
It commutes with every \(\sigma^{(\beta)}\) and with \(\mathbb E_0\).
If \(\varphi\) is a \(p\)-adic KMS form, then so is
\({}^g\varphi:=\varphi\circ\theta_g^{-1}\), and
\({}^g\varphi(1)=\varphi(1)\).
\end{proposition}

\begin{proof}
\noindent At finite level multiplication by \(g\) commutes with every
ideal action. Thus \(\theta_g\) commutes with Frobenius and preserves
the maximal \(\Lambda_K\)-order \(R_{\mathfrak f}\).
Since \(g(\mathfrak a\DR_K)=\mathfrak a\DR_K\), it also commutes
with \(\rho_{\mathfrak a}\) and hence with
\(V_{\mathfrak a,\xi}=\xi\rho_{\mathfrak a}\).
Commutation with \(\rho\) extends the action to
the rational crossed product by fixing \(U_{\mathfrak a}\) and
\(U_{\mathfrak a}^*\).  It preserves the integral generators and
hence restricts to \(A_K^{\mathrm{int}}\).
It preserves every ideal degree, hence commutes
with each \(\sigma^{(\beta)}\) and with the degree projection.
Apply the KMS condition to \(\theta_g^{-1}(X)\) and
\(\theta_g^{-1}(Y)\) to prove the last assertion.
\end{proof}

\section{Totally real fields and Deligne-Ribet \texorpdfstring{\(L\)}{L}-functions}
\label{sec:DR-padic}

\noindent For a totally real field \(K\ne\Q\), Deligne and Ribet
construct measures whose moments are differences of special values of
abelian \(L\)-functions \cite[Theorem~(8.2)]{DR80} giving rise to their theory of $p$-adic $L$-functions.
After scalar extension to \(\C_p\), the Verschiebung operators on the
periodic Witt vectors \(R_K\) act on \(E_{K,p}\). From the Deligne-Ribet
measures we construct a linear form \(\ell_\beta\) on \(E_{K,p}\) and
show that it satisfies \eqref{eq:distribution-required} for these
operators. By Proposition~\ref{prop:KMS-criterion}, this is equivalent
to the \(p\)-adic KMS condition for \(\ell_\beta\circ\mathbb E_0\)
on \(A_{K,p}\).

\subsection{Construction of $p$-adic KMS forms via Deligne-Ribet measures}
\label{subsec:DR-main}

\noindent Fix \(K\ne\Q\) totally real and \(\jmath:\Kbar\into\C_p\).
We seek a form whose value on a totally even ray character is its \(p\)-adic
\(L\)-function and whose value on \(1\) is the \(p\)-adic zeta function.
We first construct a linear form \(\ell_\beta\) on \(E_{K,p}\), then
extend it by zero on all nonidentity ideal degrees. This gives
\(\varphi_\beta=\ell_\beta\circ\mathbb E_0\), where
\(\mathbb E_0:A_{K,p}\to E_{K,p}\) is the projection
in \eqref{eq:projection-abelian-part}. By Proposition~\ref{prop:KMS-criterion},
the KMS condition for this extension amounts to
\[
 \ell_\beta(\rho_{\mathfrak a}x)
 =\bigl((\Nm(\mathfrak a))^{(\beta)}\bigr)^{-1}\ell_\beta(x)
 \qquad(\mathfrak a\in I_K^{(p)}).
\]
The Deligne-Ribet measures will give \(\ell_\beta\) by interpolation.
The factor in this formula will come from multiplying ideals by
\(\mathfrak a\) in the defining \(L\)-series.

\noindent To describe the functions to which interpolation applies,
observe that each finite \'etale algebra
\(E_{\mathfrak f}\) splits over \(\C_p\). Passing to the direct limit gives
\begin{equation}
 \Phi:E_{K,p}=\C_p\tensor_{K,\jmath}E_K
 \xrightarrow{\ \sim\ }\Clc(\DR_K,\C_p),\qquad
 \Phi(\lambda\tensor x)(\eta)=\lambda\jmath(\eta(x)).
\label{eq:padic-coefficient-splitting}
\end{equation}
Here \(\Clc\) means locally constant: every such function factors
through some \(\DR_{\mathfrak f}\). For a finite-order totally even
ray character \(\chi\) of conductor \(\mathfrak f_\chi\), let
\(x_\chi\in E_{K,p}\) be the pullback of \(\chi\) on
\(\DR_{\mathfrak f_\chi}^{\times}\), extended by zero on the
other classes at that level. In particular \(x_{\mathbf1}=1\).

\noindent Our convention for the Deligne-Ribet function
\(L_{K,p}(\beta,\chi)\) is specified by its values at the
\emph{interpolation points} \(\beta=1-k\), \(k=1,2,\ldots\):
\begin{equation}
 L_{K,p}(1-k,\chi)
 =L^{(\mathfrak f_\chi p)}(1-k,\chi\omega_p^{-k}).
\label{eq:DR-character-interpolation}
\end{equation}
Here \(\omega_p(\mathfrak a)=\omega(\Nm(\mathfrak a))\) for
\(\mathfrak a\) prime to \(p\), with \(\omega\) as in
Section~\ref{subsec:padic-dynamics}. The superscript omits Euler
factors at all primes dividing \(\mathfrak f_\chi p\Ok\), even if
the twist has smaller conductor. Character values and algebraic special values are
viewed in \(\C_p\) via \(\jmath\).

\begin{theorem}
\label{thm:DR-main-KMS}
\noindent The Deligne-Ribet measures define a family of \(p\)-adic KMS forms
\begin{equation}
 \varphi_\beta=\ell_\beta\circ\mathbb E_0:A_{K,p}\longrightarrow\C_p,
 \qquad \beta\in D_{K,p}\setminus\{1\}.
\label{eq:DR-main-linear-form}
\end{equation}
For each \(A\in A_{K,p}\), the function \(\beta\mapsto\varphi_\beta(A)\)
is meromorphic on \(D_{K,p}\), with no pole except possibly a simple
pole at \(1\). On characters and on the unit one has
\begin{equation}
 \varphi_\beta(x_\chi)=L_{K,p}(\beta,\chi),\qquad
 \varphi_\beta(1)=\zeta_{K,p}(\beta).
\label{eq:DR-partition-is-zeta}
\end{equation}
These formulas extend the usual functions from \(D_{\Q,p}\) to \(D_{K,p}\).
Where \(\zeta_{K,p}(\beta)\ne0\), the normalized form
\(\widehat\varphi_\beta=\varphi_\beta/\zeta_{K,p}(\beta)\)
is a \(p\)-adic KMS state.
\end{theorem}

\noindent At the interpolation points there is a useful nonvanishing
case. If \(k\ge2\) is divisible by \(\varphi(q_p)\), where
\(\varphi\) is Euler's totient, the finite character \(\omega_p^k\)
is trivial. Formula~\eqref{eq:DR-character-interpolation} then gives
\begin{equation}
 \zeta_{K,p}(1-k)=\zeta_K(1-k)
 \prod_{\mathfrak q\mid p}(1-\Nm(\mathfrak q)^{k-1})\ne0.
\label{eq:DR-zeta-partition}
\end{equation}
Indeed \(k\) is even, so the functional equation for the totally real
field \(K\) gives \(\zeta_K(1-k)\ne0\), and none of the displayed
Euler factors vanishes. Thus the theorem gives a state at each
of these points.

\subsection{The Deligne-Ribet measure}

\noindent We enlarge \(\DR_K\) so that multiplication by every
fractional ideal becomes invertible.  Define the localization
\[
 \mathscr I_K:=I_K^{-1}\DR_K.
\]
Its elements are formal fractions \(x/\mathfrak d\), with
\(x\in \DR_K\) and \(\mathfrak d\in I_K\), where
\[
 \frac{x}{\mathfrak d}=\frac{y}{\mathfrak e}
 \quad\Longleftrightarrow\quad
 \mathfrak e x=\mathfrak d y.
\]
Thus \(\mathfrak d^{-1}\DR_K\) denotes
\(\{x/\mathfrak d:x\in \DR_K\}\), not an inverse operation already
defined on \(\DR_K\).  Give each such copy its profinite topology
and declare these copies open.  Then
\[
 \mathscr I_K=\bigcup_{\mathfrak d\in I_K}\mathfrak d^{-1}\DR_K
\]
is locally compact, each displayed subset is compact open,
and \(J_K\) acts by homeomorphisms.

\noindent To identify this localization with the space used by Deligne and Ribet,
extend
\eqref{eq:DR-equivalence} to fractional ideals by
\[
 \mathfrak a\sim_{\mathfrak f}\mathfrak b
 \quad\Longleftrightarrow\quad
 \mathfrak a=t\mathfrak b
 \quad\text{for some }\;
 t\in K_+^\times\cap(1+\mathfrak f\mathfrak b^{-1}).
\]
Let \(\mathscr I_{\mathfrak f}=J_K/{\sim_{\mathfrak f}}\)
be the resulting discrete set, generally infinite.  Then
\[
 \mathscr I_K\cong
 \varprojlim_{\mathfrak f}\mathscr I_{\mathfrak f},
 \qquad
 \DR_K\cong\varprojlim_{\mathfrak f}\DR_{\mathfrak f}.
\]
Indeed,
\(\mathfrak a\sim_{\mathfrak f}\mathfrak b\) if and only if
\(\mathfrak d\mathfrak a\sim_{\mathfrak d\mathfrak f}
\mathfrak d\mathfrak b\), and every compatible family has a
common denominator: a denominator for a representative of its component at
\((1)\) works for every component
\cite[(2.2), (2.5), (2.8) and p.~240]{DR80}.
Our \(\DR_K\) and \(\mathscr I_K\) are respectively their spaces
\(A\) and \(I\).

\noindent The norm of
\cite[(2.12)]{DR80} has a continuous \(p\)-component
\(\mathcal N_p:\mathscr I_K\to\Q_p\), with
\begin{equation}
 \mathcal N_p(\DR_K)\subseteq\Z_p,\qquad
 \mathcal N_p(\mathfrak a x)=\Nm(\mathfrak a)\mathcal N_p(x).
\label{eq:DR-norm}
\end{equation}
Its restriction to \(\DR_K^\times\) takes values in
\(\Z_p^\times\) and has open image.
Only the following compact open set will enter the integral:
\begin{equation}
 U_p:=\{x\in \DR_K:\mathcal N_p(x)\in\Z_p^\times\}.
\label{eq:U-p}
\end{equation}
It is stable under \(I_K^{(p)}\) and contains noninvertible points
away from \(p\).

\begin{lemma}
\label{lem:norm-image}
\noindent The norm image of the integration set is
\[
 \mathcal N_p(U_p)=H_{K,p}.
\]
In particular,
\(\left|\log_p\langle\mathcal N_p(y)\rangle\right|_p
\le p^{-m_{K,p}}\) for every \(y\in U_p\).
\end{lemma}

\begin{proof}
\noindent Integral ideals are dense in \(\DR_K\), and \(U_p\)
is clopen with \(I_K\cap U_p=I_K^{(p)}\).
Thus prime-to-\(p\) integral ideals are dense in \(U_p\).
Continuity and compactness give
\(\mathcal N_p(U_p)=\overline{\Nm(I_K^{(p)})}\).
A closed submonoid of a profinite group is a subgroup:
positive powers of each element approximate its inverse in
the finite quotients.  This closure is therefore
\(\overline{\Nm(J_K^{(p)})}=H_{K,p}\).
The bound follows from \eqref{eq:norm-subgroup}.
\end{proof}

\begin{remark}
\label{rem:measure-support}
\noindent The distinction between \(U_p\) and \(\DR_K^\times\) is essential.
For a proper \(\mathfrak a\in I_K^{(p)}\), the image
\(\mathfrak a\DR_K\) misses every unit. A form supported only on
units would vanish on every \(\rho_{\mathfrak a}x\), contrary to
\eqref{eq:distribution-required} for a nonzero form.
In contrast, if \(\mathfrak q\mid p\), the idempotent
\(e_{\mathfrak q}\) vanishes on \(U_p\). Its value under any form
constructed by integration there is zero, explaining why we do
not impose covariance for primes above \(p\).
\end{remark}

\noindent For a locally constant, compactly supported test function
\(\varepsilon\) on \(\mathscr I_K\) with algebraic values, let
\begin{equation}
 L(s,\varepsilon)=\sum_{\mathfrak b\in J_K}
       \varepsilon(\mathfrak b)\Nm(\mathfrak b)^{-s}
       \qquad(\Re(s)>1).
\label{eq:DR-L-function}
\end{equation}
We use a complex embedding for this series. Its special values
\(L(1-k,\varepsilon)\), \(k\ge1\), are algebraic by
\cite[Theorem~(8.1)]{DR80}. For \(\C_p\)-valued test functions
use \(\C_p\)-linear extension from characteristic functions.
For each unit \(c\in\DR_K^\times\), Deligne and Ribet construct
a bounded \(p\)-adic measure \(\mu_c\) (the \(p\)-component of their
integral measure) whose moments are
\begin{equation}
 \int_{\mathscr I_K}\varepsilon(y)\mathcal N_p(y)^{k-1}\,d\mu_c(y)
 =L(1-k,\varepsilon)
 -\mathcal N_p(c)^k L(1-k,\varepsilon(c\,\cdot)).
\label{eq:DR-moment}
\end{equation}
This is \cite[Theorem~(8.2)]{DR80}. The second term is the
correction that makes the distribution a measure. We remove it
one finite level at a time.

\subsection{$p$-adic KMS forms and proof of Theorem~\ref{thm:DR-main-KMS}}

\noindent The measure \(\mu_c\) in \eqref{eq:DR-moment} gives a
\emph{difference} of special values. We first recover the individual
values. The point is that each locally constant test function factors
through a finite ray level. We may choose the auxiliary unit \(c\)
to act trivially at that level. Smoothing then multiplies the special
value by a scalar, which we can divide out.

\noindent Extend functions on \(U_p\) by zero to \(\mathscr I_K\).
Put \(\omega_p=\omega\circ\mathcal N_p\) on \(U_p\) and
\(w_\beta(u)=(u^{(\beta)})^{-1}\) for \(u\in H_{K,p}\).
For \(c\in\DR_K^\times\), take the Mellin transform of the measure
\(f\,\mu_c\) along the norm \cite[p.~178 and \S4]{Ribet79}:
\[
 M_{\beta,c}(f)=\int_{U_p} f(y)\,w_\beta(\mathcal N_p(y))\,d\mu_c(y).
\]
This means integrating the character \(w_\beta\) against the image
of \(f\,\mu_c\) under \(\mathcal N_p\).

\begin{lemma}
\label{lem:DR-test-interpolation}
\noindent For each locally constant \(f:U_p\to\C_p\) there is a
unique rigid meromorphic function \(\ell_\beta(f)\) on \(D_{K,p}\)
with at most a simple pole at \(1\), no other poles, and values
\begin{equation}
 \ell_{1-k}(f)=L(1-k,f\omega_p^{-k})\qquad(k\ge1).
\label{eq:DR-interpolation}
\end{equation}
For fixed \(\beta\ne1\), this assignment is linear in \(f\).
\end{lemma}

\begin{proof}
\noindent The extension of \(f\) to \(\DR_K\) factors through some
\(\DR_{\mathfrak f}\). Choose a unit \(c\) acting trivially at that
level and with \(\log_p\mathcal N_p(c)\ne0\). Such a unit exists:
the cyclotomic norm on \(\DR_K^\times\) has open image, and raising a
unit of nonzero logarithm to the order of its image at level
\(\mathfrak f\) does not kill its logarithm. Thus \(f(cy)=f(y)\).

\noindent Write \(u_c=\mathcal N_p(c)\). At \(\beta=1-k\) we have
\(w_{1-k}(u)=u^{k-1}\omega(u)^{-k}\).
Apply \eqref{eq:DR-moment} to \(\varepsilon=f\omega_p^{-k}\).
Since \(f(cy)=f(y)\), we have
\(\varepsilon(cy)=\omega(u_c)^{-k}\varepsilon(y)\), and hence
\[
 M_{1-k,c}(f)
   =(1-\langle u_c\rangle^k)L(1-k,f\omega_p^{-k}).
\]
This identifies the scalar to remove. For \(\beta\ne1\) put
\begin{equation}
 \ell_\beta(f)=\frac{M_{\beta,c}(f)}
      {1-\exp\bigl((1-\beta)\log_p u_c\bigr)} .
\label{eq:DR-measure-linear-form}
\end{equation}
Lemma~\ref{lem:norm-image} gives uniform convergence of the
exponential series in the integrand on each closed subdisc of
\(D_{K,p}\), so the numerator is analytic there.
Since \(1-\beta\in D_{K,p}\), the same bound puts
\((1-\beta)\log_p u_c\) in the disc of injectivity of the exponential.
The denominator therefore vanishes only at \(1\), with nonzero
derivative \(\log_p u_c\). At \(1-k\) it is \(1-\langle u_c\rangle^k\),
which proves \eqref{eq:DR-interpolation}.

\noindent Although the choice of \(c\) may depend on the finite level,
the quotient does not depend on the choice.
The interpolation points include \(1-(1+p^n)=-p^n\), which
accumulate at \(0\). The identity principle on the disc proves
uniqueness (apply it after multiplying by \(\beta-1\)).
It gives independence of \(c\) and linearity, since the
interpolation values have these properties. It is the smoothed
measure \(\mu_c|_{U_p}\) that is bounded. After division we need
only a linear form on locally constant functions.
\end{proof}

\begin{proof}[Proof of Theorem~\ref{thm:DR-main-KMS}]
\noindent For \(x\in E_{K,p}\), apply the lemma to
\(f_x=\Phi(x)|_{U_p}\) and write \(\ell_\beta(x)=\ell_\beta(f_x)\).
It remains to check the KMS condition and identify the character values.
For \(\mathfrak a\in I_K^{(p)}\), formula~\eqref{eq:rho-functions}
gives \(f_{\rho_{\mathfrak a}x}(\mathfrak a y)=f_x(y)\) on \(U_p\),
and \(f_{\rho_{\mathfrak a}x}\) vanishes outside \(\mathfrak a U_p\).
First take \(f_x\) with algebraic values. In the Dirichlet series
\eqref{eq:DR-L-function}, only ideals
\(\mathfrak a\mathfrak b\) contribute, with \(\mathfrak b\) integral
and prime to \(p\). Changing the index to \(\mathfrak b\) gives
\[
 L(s,f_{\rho_{\mathfrak a}x}\omega_p^{-k})
 =\Nm(\mathfrak a)^{-s}\omega_p(\mathfrak a)^{-k}
    L(s,f_x\omega_p^{-k}).
\]
The identity continues to \(s=1-k\), where the multiplier is
\(w_{1-k}(\Nm(\mathfrak a))\). At these special values the identity
holds for arbitrary \(f_x\) by \(\C_p\)-linear extension from
characteristic functions. Thus the two sides of the following
formula agree at all \(1-k\), and the lemma and the identity
principle give it throughout \(D_{K,p}\setminus\{1\}\):
\begin{equation}
 \ell_\beta(\rho_{\mathfrak a}x)=w_\beta(\Nm(\mathfrak a))\ell_\beta(x).
\label{eq:DR-rho-covariance}
\end{equation}
Since \(V_{\mathfrak a,\xi}=\jmath(\xi)\rho_{\mathfrak a}\) after
scalar extension, this is the covariance criterion
\eqref{eq:distribution-required}. Hence \(\ell_\beta\circ\mathbb E_0\)
satisfies the \(p\)-adic KMS condition.

\noindent Finally \(f_{x_\chi}\) on ideal points is \(\chi\) on
ideals prime to \(\mathfrak f_\chi p\) and zero elsewhere.
Thus \eqref{eq:DR-interpolation} gives exactly
\eqref{eq:DR-character-interpolation}. Uniqueness on \(D_{\Q,p}\)
identifies these values with the usual Deligne-Ribet function,
and the lemma continues it to \(D_{K,p}\). The trivial character
gives the value on \(1\). Dividing by this value when it is nonzero
preserves the KMS condition and gives a state.
\end{proof}

\noindent On a ray character one can also see the smoothing factor
for an \emph{arbitrary} auxiliary unit \(c\). Put
\(f_\chi=\Phi(x_\chi)|_{U_p}\) and \(u_c=\mathcal N_p(c)\).
Since \(f_\chi(cy)=\chi(c)f_\chi(y)\), the moment formula at
\(1-k\) and the identity principle give the meromorphic identity
\[
 M_{\beta,c}(f_\chi)
 =\bigl(1-\chi(c)\exp((1-\beta)\log_p u_c)\bigr)
       L_{K,p}(\beta,\chi).
\]
For the unit chosen in the lemma, \(\chi(c)=1\) and
\(\log_p u_c\ne0\). Thus for \(\beta\ne1\) we can divide by this
factor to recover \(L_{K,p}\).
Dividing afterwards by \(\zeta_{K,p}\), where nonzero and regular,
serves a different purpose: it makes the KMS form take the value \(1\)
on the unit.

\subsection{KMS states and Galois symmetries}

\noindent The normal form \eqref{eq:padic-normal-form} reduces the
evaluation of a state constructed above to its restriction to
\(E_{K,p}\). The first formula below makes this reduction explicit.
The second describes how translating the state changes its character
values.

\filbreak
\begin{proposition}
\label{prop:state-values}
\noindent Fix \(\beta\in D_{K,p}\setminus\{1\}\) with
\(\zeta_{K,p}(\beta)\ne0\).
\begin{enumerate}[label=(\roman*),leftmargin=*,itemsep=3pt]
\item For \(\mathfrak a,\mathfrak b\in I_K^{(p)}\),
\(\xi\in\mathfrak a\), and \(x\in E_{K,p}\),
\begin{equation}
 \widehat\varphi_\beta((\xi U_{\mathfrak a})xU_{\mathfrak b}^*)
 =
 \begin{cases}
 \displaystyle
 \frac{\jmath(\xi)}{(\Nm(\mathfrak a))^{(\beta)}}
       \widehat\varphi_\beta(x),&\mathfrak a=\mathfrak b,\\[4pt]
 0,&\mathfrak a\ne\mathfrak b.
 \end{cases}
\label{eq:state-monomials}
\end{equation}
The indices need not be coprime.  In particular,
\[
 \widehat\varphi_\beta(e_{\mathfrak a})
      =((\Nm(\mathfrak a))^{(\beta)})^{-1},
 \qquad e_{\mathfrak a}=U_{\mathfrak a}U_{\mathfrak a}^*.
\]
\item For \(g\in G_K^{\ab}\) and a finite-order totally even ray
character \(\chi\), the translated state
\({}^g\widehat\varphi_\beta
=\widehat\varphi_\beta\circ\theta_g^{-1}\) satisfies
\begin{equation}
 {}^g\widehat\varphi_\beta(x_\chi)
  =\chi(g)^{-1}
       \frac{L_{K,p}(\beta,\chi)}{\zeta_{K,p}(\beta)}.
\label{eq:state-Galois-character}
\end{equation}
\end{enumerate}
\end{proposition}

\begin{proof}
\noindent If \(\mathfrak a\ne\mathfrak b\), the monomial has
nonidentity ideal degree and is killed by \(\mathbb E_0\).
If \(\mathfrak a=\mathfrak b\), it equals
\(V_{\mathfrak a,\xi}(x)\), so
\eqref{eq:distribution-required} gives the first formula.
Taking \(x=1\) and dividing by any nonzero \(\jmath(\xi)\)
gives the idempotent value.
Finally, \(\theta_gx_\chi=\chi(g)x_\chi\). Apply the definition of
the action on forms and Theorem~\ref{thm:DR-main-KMS} to obtain (ii).
\end{proof}

\begin{remark}
\noindent The table compares the classical system
\cite[Theorem~2.1]{LLN09} with the family constructed here. In the
classical column \(\psi_{\beta,w}\) denotes the extremal state
indexed by \(w\in\DR_K^\times\) at real \(\beta>1\).
In the \(p\)-adic column \(\widehat\varphi_\beta\) is defined wherever
the form is regular and its value on \(1\) is nonzero. \\ 
It would be interesting to understand conceptually the precise interaction of our $p$-adic BC-systems with (local) class field theory.
\begin{center}
{\small\setlength{\tabcolsep}{4pt}
\begin{tabular}{@{}>{\raggedright\arraybackslash}p{0.20\textwidth}>{\raggedright\arraybackslash}p{0.34\textwidth}>{\raggedright\arraybackslash}p{0.39\textwidth}@{}}
\toprule
& Classical BC-system & \(p\)-adic BC-system\\
\midrule
Partition function &
\(\zeta_K(\beta)\) &
\(\zeta_{K,p}(\beta)=\varphi_\beta(1)\)\\[4pt]
Galois symmetries &
\(G_K^{\ab}\) commutes with \(\sigma_t\). &
\(G_K^{\ab}\) commutes with each \(\sigma^{(\beta)}\).\\[4pt]
KMS states and forms &
Positive normalized complex KMS states &
\(\C_p\)-linear KMS forms, normalized by division by \(\zeta_{K,p}(\beta)\)\\[4pt]
Galois compatibility &
\({}^g\psi_{\beta,w}=\psi_{\beta,g^{-1}w}\) &
\({}^g\widehat\varphi_\beta(x_\chi)=
 \chi(g)^{-1}\widehat\varphi_\beta(x_\chi)\)\\
\bottomrule
\end{tabular}}
\end{center}
\noindent In both columns \({}^g\psi=\psi\circ\theta_g^{-1}\):
the Galois group acts on the argument of the state, not on its scalars.
\end{remark}

\subsection{Real quadratic fields}

\noindent Let \(K\) be real quadratic of discriminant \(\Delta\),
and let \(\chi_\Delta\) be the associated primitive quadratic Dirichlet
character. At the integers below there is no finite-character twist,
so the classical factorization \(\zeta_K(s)=\zeta(s)L(s,\chi_\Delta)\)
\cite[p.~1150, (6)]{Zagier99} applies directly to
\eqref{eq:DR-zeta-partition}.

\begin{remark}
\label{prop:quadratic-partition}
\noindent If \(k\ge2\) is divisible by \(\varphi(q_p)\), then
\begin{equation}
 \zeta_{K,p}(1-k)
 =\zeta(1-k)L(1-k,\chi_\Delta)
       \prod_{\mathfrak q\mid p}
          \bigl(1-\Nm(\mathfrak q)^{k-1}\bigr).
\label{eq:quadratic-partition}
\end{equation}
For each prime \(\mathfrak q\mid p\), its residue degree is
\(f_{\mathfrak q}=[\Ok/\mathfrak q:\F_p]\), so
\(\Nm(\mathfrak q)=p^{f_{\mathfrak q}}\).
The factorization of \(p\Ok\) gives:
\begin{center}
\begin{tabular}{@{}llll@{}}
\toprule
Type & \(p\Ok\) & Residue degrees & Euler multiplier\\
\midrule
Split &
\(\mathfrak q\mathfrak q'\), \(\mathfrak q\ne\mathfrak q'\) &
\(f_{\mathfrak q}=f_{\mathfrak q'}=1\) & \((1-p^{k-1})^2\)\\[3pt]
Inert &
\(\mathfrak q\) &
\(f_{\mathfrak q}=2\) & \(1-p^{2(k-1)}\)\\[3pt]
Ramified &
\(\mathfrak q^2\) &
\(f_{\mathfrak q}=1\) & \(1-p^{k-1}\)\\
\bottomrule
\end{tabular}
\end{center}
There is one factor per distinct prime above \(p\), irrespective of
its ramification index. The table includes \(p=2\), when the condition
on \(k\) says simply that \(k\) is even.
\end{remark}

\begin{corollary}
\label{cor:quadratic-weight-one}
\noindent The normalized family \(\widehat\varphi_\beta\)
extends pointwise analytically to \(\beta=1\), where it is a
\(p\)-adic KMS state.  Its value on \(A\in A_{K,p}\) is
\[
 \widehat\varphi_1(A)=
 \frac{\operatorname*{Res}_{\beta=1}\varphi_\beta(A)}
      {\operatorname*{Res}_{\beta=1}\zeta_{K,p}(\beta)}.
\]
In particular,
\(\widehat\varphi_1(e_{\mathfrak a})=(\Nm(\mathfrak a))^{-1}\)
for \(\mathfrak a\in I_K^{(p)}\).
\end{corollary}

\begin{proof}
\noindent Leopoldt's assertion holds in the rank-one unit case.
Indeed, a suitable power of a fundamental unit lies in a
principal-unit subgroup on which the \(p\)-adic logarithm is
injective. Its logarithm is nonzero because the unit is not torsion.  The \(p\)-adic regulator is therefore nonzero.
The residue formula of \cite{Colmez88} implies that
\(\zeta_{K,p}\) has a simple pole at \(1\) with nonzero residue.

\noindent For every fixed \(A\), the numerator
\(\varphi_\beta(A)\) has at most a simple pole by
Theorem~\ref{thm:DR-main-KMS}.
Consequently
\[
 \widehat\varphi_\beta(A)=
 \frac{(\beta-1)\varphi_\beta(A)}
      {(\beta-1)\zeta_{K,p}(\beta)}
\]
is analytic near \(1\), and its value there is the stated ratio
of residues.  Its value on the unit is \(1\).
For homogeneous \(Y\), the automorphism \(\sigma^{(\beta)}\)
acts by an analytic scalar. Taking residues in the KMS condition
therefore proves it at \(1\) for the residue form and its normalization.
Linearity gives the condition for arbitrary \(Y\).  Finally, the formula for \(e_{\mathfrak a}\)
follows by continuity from \eqref{eq:state-monomials}, since
\(u^{(1)}=u\).
\end{proof}

\subsection{$p$-adic KMS states for CM number fields}

\noindent 
For a CM field equipped with a \(p\)-ordinary CM type, Katz constructed \(p\)-adic measures and \(p\)-adic \(L\)-functions
\cite{Katz78,HT93}. The corresponding BC-systems have been studied in \cite{Yal12CM} and \cite{UraIII,UraIV}. It's a natural question whether the methods of this paper can be extended to construct \(p\)-adic KMS states for CM fields, and how they relate to Katz's measures. This remains an open problem for further research.

\begingroup
\renewcommand{\addcontentsline}[3]{}
\section*{Use of AI tools}

\noindent We used LLMs to  improve clarity, grammar, and style of the text, as well as to check for consistency in notation and formatting. The LLM did not contribute to the mathematical content or results of the paper.
Responsibility for the final mathematical claims and attributions
remains with the author.
\par
\endgroup

\bibliographystyle{alpha}
\bibliography{int_padic_bc}

\end{document}